\documentclass[12pt]{article}
\usepackage{amsmath,amsfonts,amssymb,amscd}
\usepackage{graphicx}
\usepackage[T2A]{fontenc}

\newtheorem{lem}{Lemma}[section]
\newtheorem{prop}{Proposition}[section]

\newtheorem{cor}{Corollary}[section]
\newtheorem{rem}{Remark}[section]
\newtheorem{exa}{Example}[section]
\newtheorem{dfn}{Definition}[section]

\makeatletter \@addtoreset{equation}{section} \makeatother

\newcommand{\mC}{\mathbb{C}}

\newcommand{\mR}{\mathbb{R}}
\newcommand{\mT}{\mathbb{T}}
\newcommand{\mZ}{\mathbb{Z}}
\newcommand{\mN}{\mathbb{N}}

\newcommand{\one}{{\bf 1}}

\newcommand{\bX}{{\bf X}}

\newcommand{\bc}{{\bf c}}
\newcommand{\bd}{{\bf d}}

\newcommand{\bT}{{\bf T}}

\newcommand{\calB}{{\cal B}}

\newcommand{\calH}{{\cal H}}

\newcommand{\calM}{{\cal M}}
\newcommand{\calN}{{\cal N}}
\newcommand{\calO}{{\cal O}}

\newcommand{\calX}{{\cal X}}

\newcommand{\DT}{{\cal DT}}

\newcommand{\frakh}{{\mathfrak h}}
\newcommand{\frakX}{{\mathfrak X}}

\newcommand{\eps}{\varepsilon}
\newcommand{\ph}{\varphi}
\newcommand{\thet}{\vartheta}

\newcommand{\auto}{\operatorname{auto}}
\newcommand{\Conv}{\operatorname{Conv}}
\newcommand{\Koop}{\operatorname{Koop}}

\newcommand{\dist}{\operatorname{dist}}
\newcommand{\diag}{\operatorname{diag}}
\newcommand{\diam}{\operatorname{diam}}

\newcommand{\id}{\operatorname{id}}

\newcommand\qed{{\unskip\nobreak\hfil\penalty50
  \hskip2em\hbox{}\nobreak\hfil\mbox{\rule{1ex}{1ex} \qquad}
    \parfillskip=0pt \finalhyphendemerits=0\par\medskip}}

\begin{document}

\title
{$\mu$-norm of an operator}
\author{D.Treschev \\
Steklov Mathematical Institute of Russian Academy of Sciences
}
\date{}
\maketitle

\begin{abstract}
Let $(\calX,\mu)$ be a measure space. For any measurable set $Y\subset\calX$ let $\one_Y : \calX\to\mR$ be the indicator of $Y$ and let $\pi_Y$ be the orthogonal projector $L^2(\calX)\ni f\mapsto\pi_Y f = \one_Y f$. For any bounded operator $W$ on $L^2(\calX,\mu)$ we define its $\mu$-norm
$\|W\|_\mu = \inf_\chi\sqrt{\sum \mu(Y_j) \|W\pi_Y\|^2}$, where the infinum is taken over all measurable partitions $\chi = \{Y_1,\ldots,Y_J\}$ of $\calX$.  We present some properties of the $\mu$-norm and some computations. Our main motivation is the problem of the construction of a quantum entropy.
\end{abstract}

\section{Introduction and motivation}
\label{sec:intro}

Let $\calX$ be a nonempty set and let $\calB$ be a $\sigma$-algebra of subsets $X\subset\calX$. Consider the measure space $(\calX,\calB,\mu)$, where $\mu$ is a probability measure: $\mu(\calX) = 1$.

Consider the Hilbert space $\calH = L^2(\calX,\mu)$ with the scalar product and the norm
$$
  \langle f,g\rangle = \int_\calX f\overline g \, d\mu, \quad
  \|f\| = \sqrt{\langle f,f\rangle} .
$$
For any bounded operator $W$ on $\calH$ let $\|W\|$ be its norm defined by
$$
  \|W\| = \sup_{\|f\|=1}  \|Wf\|.
$$

For any $X\in\calB$ consider the orthogonal projector
\begin{equation}
\label{imath}
  \pi_X : \calH\to\calH, \qquad
  \calH\ni f \mapsto \pi_X f = \one_X \cdot f,
\end{equation}
where $\one_X$ is the indicator of $X$ (the function equal to 1 on $X$ and vanishing at other points). Then $\mu(X)>0$ implies $\|\pi_X\| = 1$ and for any $X',X''\in\calB$
$$
  \pi_{X'} + \pi_{X''} = \pi_{X'\cup X''} + \pi_{X'\cap X''}, \quad
  \pi_{X'} \pi_{X''} = \pi_{X'\cap X''}.
$$
We say that $\chi = \{Y_1,\ldots,Y_J\}$ is a (finite, measurable) partition (of $\calX$) if
$$
  Y_j\in\calB, \quad
  \mu\big(\calX\setminus \cup_{1\le j\le J} Y_j\big) = 0, \quad
  \mu(Y_j\cap Y_k) = 0 \quad
  \mbox{ for any $j,k\in\{1,\ldots,J\}$, $k\ne j$}.
$$

For any two partitions $\chi = \{Y_1,\ldots,Y_J\}$ and $\kappa = \{X_1,\ldots,X_K\}$ we denote
$$
  \chi\vee\kappa = \{Y_j\cap X_k\}_{j=1,\ldots,J,\, k=1,\ldots,K}.
$$
Obviously, $\chi\vee\kappa$ is also a partition.

Let $W$ be a bounded operator on $\calH$. For any partition $\chi = \{Y_1,\ldots,Y_J\}$ of $\calX$ we define
\begin{equation}
\label{calM}
  \calM_\chi(W) = \sum_{j=1}^J \mu(Y_j) \| W\pi_{Y_j} \|^2.
\end{equation}
We define the $\mu$-norm\footnote
{in fact, a seminorm}
of $W$ by
\begin{equation}
\label{mumeasure}
  \|W\|_\mu = \inf_\chi \sqrt{\calM_\chi(W)}.
\end{equation}

Recall that the operator $U$ is said to be an isometry if
$$
  \langle f,g\rangle = \langle Uf,Ug\rangle, \qquad
  f,g\in\calH.
$$
If the isometry $U$ is invertible then $U$ is called a unitary operator.

For any bounded $W$, any $Y\in\calB$, and any isometry $U$
$$
  \|W\pi_Y\| \le \|W\|, \quad
  \|UW\| = \|W\|.
$$
This implies the following obvious properties of the $\mu$-norm:
\begin{eqnarray}
\label{|1|}
         \|\id\|_\mu
     &=& 1,\qquad
         \|W\|_\mu
 \;\le\; \|W\|,     \\
\label{WW}
          \|W_1 W_2\|_\mu
    &\le& \|W_1\| \|W_2\|_\mu, \\
\label{lambdaW}
         \|\lambda W\|_\mu
     &=& |\lambda|\, \|W\|_\mu \quad
         \mbox{ for any $\lambda\in\mC$}, \\
\label{UWWU}
         \|W\|_\mu
     &=& \|UW\|_\mu \quad \mbox{ for any unitary $U$}.
\end{eqnarray}
\medskip

The first question which probably comes to the reader's mind is ``why?''. Why such a construction may be useful or interesting? Now we are going to explain our motivations.

Let $F:\calX\to\calX$ be an endomorphism of the measure space $\calX,\calB,\mu$. This means that for any $X\in\calB$ the set $F^{-1}(X)$ (the complete preimage) also lies in $\calB$ and $\mu(X) = \mu(F^{-1}(X))$. Invertible endomorphisms are called automorphisms. Let $End(\calX)$ denote the semigroup of all endomorphisms of $(\calX,\calB,\mu)$. There are two standard constructions associated with any
$F\in End(\calX,\mu)$.
\medskip

(1). Any such $F$ generates the isometry (a unitary operator if $F$ is an automotphism) $U_F$ on $L^2(\calX,\mu)$ (the Koopman operator):
$$
  L^2(\calX,\mu)\ni f \mapsto U_F f = f\circ F, \qquad
  U_F = \Koop(F).
$$

(2) For any $F\in End(\calX,\mu)$ it is possible to compute the measure entropy (another name is the Kolmogorov-Sinai entropy) $h(F)$.
\medskip

Our question is as follows. Is it possible to determine in some ``natural way'' a real nonnegative function on the semigroup $Iso(L^2(\calX,\mu))$ so that the diagram
$$
\begin{array}{rcl}
              &    End(\calX,\mu)                  &      \\
  h\swarrow   &                                    &  \searrow\Koop \\
 \mR_+        & \stackrel{\frakh}{\longleftarrow}  &  Iso(L^2(\calX,\mu))
\end{array}
$$
is commutative.

Recall the construction of the measure entropy of an endomorphism. Let $J_N$ be the set of multiindices $j = (j_0,\ldots,j_N)$, where any component $j_n$ takes values in the set $\{0,\ldots,K\}$. For any partition $\chi = \{X_0,\ldots,X_K\}$ and $j\in J_N$ we define
$$
  \bX_j = F^{-N}(X_{j_N})\cap \ldots \cap F^{-1}(X_{j_1}) \cap X_{j_0}.
$$
We define $h_F(\chi,N+1)$ by
$$
    h_F(\chi,N+1)
  = - \sum_{j\in J_N} \mu(\bX_j) \log \mu(\bX_j).
$$
The function $h_F$, as a function of the second argument, is subadditive:
$h_F(\chi,n+m) \le h_F(\chi,n) + h_F(\chi,m)$. This implies existence of the limit
$$
  h_F(\chi) = \lim_{n\to\infty} \frac1n h_F(\chi,n).
$$
Finally, the measure entropy is defined by
$$
  h(F) = \sup_\chi h_F(\chi).
$$

Our idea is to construct the entropy of a unitary operator $U$ analogously with the following difference. Instead of $\bX_j$ we take
$$
    \frakX_j
  = \pi_{X_{j_N}} U \pi_{X_{j_{N-1}}} U \ldots U \pi_{X_{j_1}} U \pi_{X_{j_0}} .
$$
We define
$$
    \frakh_U(\chi,N+1)
  = - \sum_{j\in J_N} \|\frakX_j\|_\mu^2 \log \|\frakX_j\|_\mu^2 .
$$
Other details are the same:
\begin{equation}
\label{limfrakh}
    \frakh_U(\chi)
  = \lim_{n\to\infty} \frac1n \frakh(\chi,n), \qquad
    \frakh(U)
  = \sup_\chi \frakh_U(\chi).
\end{equation}

Corollary \ref{cor:piUpiU} (below) implies that for any automorphism $F$
\begin{equation}
\label{UF=F}
  \frakh(U_F) = h(F) .
\end{equation}

It is interesting to see what this construction gives in finite-dimensional case. Suppose $\calX$ is finite, $\#\calX = J$. Then the space $\calH$ is isomorphic to $\mC^J$. Operators on $\calH$ are identified with $J\times J$ matrices. By using results of Section \ref{subsec:finite}, it is possible to compute for any unitary
$U = (U_{jk})$
\begin{equation}
\label{h(finite)}
      \frakh(U)
   =  -\frac1J \sum_{\sigma_0,\sigma_1=1}^J
       |U_{\sigma_0\sigma_1}|^2\, \log |U_{\sigma_0\sigma_1}|^2.
\end{equation}
Note that the numbers $P_{jk} = |U_{jk}|^2$ can be regarded as elements of a bi-stochastic matrix $P$. Then the vector $\nu = (1/J,\ldots,1/J)$ is a stationary distribution for the corresponding Markov chain. It is known \cite{CT} that if a Markov chain is irreducible, aperiodic and positive recurrent then its stationary distribution is unique and its entropy rate is determined by
$-\sum_{j,k} \mu_j P_{jk} \log P_{jk}$. This quantity coincides with $\frakh(U)$.

In the literature there exist several attempts to extend the concept of the measure entropy to quantum systems, see \cite{CNT,Ohya95,Ohya00,AOW,AF,M,GLW} and many others. In \cite{A} several mutual relations between these approaches are given. In some definitions of quantum entropy (see, for example, \cite{AOW, Sr, Pe, BG, KK}) equations similar to (\ref{h(finite)}) appear. In \cite{DF05, DF18} a construction for measure entropy is proposed for doubly stochastic operators on various spaces of functions on a measure space. The question, which approach to quantum generalization of the measure entropy is ``more physical'' remains unclear.
\medskip

Concerning our definition of $\frakh(U)$, many technical questions appear, including the question on subadditivity\footnote
{This subadditivity is important for the existence of the limit (\ref{limfrakh}).}
of $\frakh_U(\chi,n)$, on the possibility to include into the construction isometries $U$ together with unitary operators, on methods of computation (or at least, on effective lower and upper estimates) of the quantities $\frakh(U)$, and many others. Answers to these questions depend on the detailed analysis of the operation $\|\cdot\|_\mu$. We start the analysis in this paper.
\medskip

The main results we present in this paper are as follows (below $W$ is an arbitrary bounded operator).

\begin{itemize}
\item $\|\pi_X\|_\mu^2 = \mu(X)$ for any $X\in\calB$ (Lemma \ref{lem:mupi}).

\item If $\chi'$ is a subpartition of $\chi$ then $\calM_{\chi'}(W) \le \calM_\chi(W)$ (Lemma \ref{lem:muM} and Corollary \ref{cor:subpartition}). Hence the quantities $\calM_\chi(W)$ approach the infinum (\ref{mumeasure}) on fine partitions.

\item For any two bounded operators $W_1$ and $W_2$ we prove that
$$
  \|W_1 + W_2\|_\mu \le \|W_1\|_\mu + \|W_2\|_\mu
$$
(Lemma \ref{lem:tri}). This triangle inequality combined with (\ref{lambdaW}) imply that $\|\cdot\|_\mu$ is a seminorm on the space of bounded operators on $\calH$.

\item Let $F$ be an automorphism of $(\calX,\calB,\mu)$ and let $U_F$ be the corresponding unitary operator on $\calH$, defined by $f\mapsto U_F = f\circ F$. Then
$$
  \|W U_F\|_\mu = \|W\|_\mu
$$
  (Lemma \ref{lem:UWU}). This implies $\|U_F^{-1} W U_F\|_\mu = \|W\|_\mu$. Informally speaking, this means that measure preserving coordinate changes on $\calX$ preserve the $\mu$-norm.

\item $\|\cdot\|_\mu$ is a continuous function in the $L^2(\calX)$ operator topology (Corollary \ref{cor:continuous}).

\item If the measure $\mu$ has no atoms then $\|W+W_0\|_\mu=\|W\|_\mu$ for any bounded $W$ and compact $W_0$ (Corollary \ref{cor:B+K}).

\item Given $g\in L^\infty(\calX)$ let $\widehat g$ be the operator defined by $f\mapsto\widehat g f = gf$. Then $\|\widehat g\|_\mu = \|g\|$, (Lemma \ref{lem:W_g}).

\item Suppose $\calX = \{1,\ldots,J\}$ is finite and the measure of any element equals $1/J$. Let $W = (W_{jk})_{j,k=1,\ldots,J}$ be an operator on $\calH \cong \mC^J$. Then by Lemma \ref{lem:finite}
$$
  \|W\|_\mu^2 = \frac1J \sum_{j,k=1}^J |W_{jk}|^2 .
$$

\item Let $\pi_H : \calH\to H$ be the orthogonal projector to the closed subspace $H\subset\calH$. We define the $\mu$-dimension of $H$ by
    $\dim_\mu(H) := \|\pi_H\|_\mu^2$. According to Lemma \ref{lem:subspace} the following statements hold.

1. If $\mu$ has no atoms and $H$ is finite-dimensional then $\dim_\mu(H)=0$.

2. If $H = \pi_X(\calH)$ for some $X\in\calB$ then $\dim_\mu(H) = \mu(X)$.

3. If $\calX$ is finite and $\mu$ is uniformly distributed among the points then $\calH = \mC^J$, $J = \#\calX$ and for any subspace $H\subset\calH\;$
$\dim_\mu(H) = \dim(H) / J$, where $\dim$ is the ordinary dimension.
\smallskip

For any $\mu$-preserving almost free action of the cyclic group $\mZ_q$ on $\calX$ we consider the space $H_0\subset\calH$ of $\mZ_q$-invariant functions. It is natural to call these functions periodic with period $1/q$. We prove that $\dim_\mu (H_0) = 1/q$ (Lemma \ref{lem:cyclic}).

\item For any partition $\{X_1,\ldots,X_K\}$
$$
  \|W\|_\mu^2 = \sum_{k=1}^K \|W \pi_{X_k}\|_\mu^2
$$
  (the right additivity of $\|\cdot\|_\mu^2$, Lemma \ref{lem:=}).

\item Instead of the left additivity we have the weaker property
  (Lemma \ref{lem:<=}):
$$
  \|W\|_\mu^2 \le \sum_{k=1}^K \|\pi_{X_k} W\|_\mu^2 .
$$

\item In Section \ref{sec:add} we consider the case when $\calX$ is a compact metric space and $\mu$ is a Borel measure w.r.t. the corresponding topology. Let $B_r(x)\subset\calX$ denote the open ball with center at $x$ and radius $r$. We prove (Proposition \ref{prop:theta}) that for any $x$ the limit
$$
  \thet(x) = \lim_{\eps\searrow 0} \|W\pi_{B_\eps(x)}\|^2
$$
  exists, the function $\thet : \calX\to\mR$ is measurable and
  $\|W\|_\mu^2 \le \int_\calX \thet \,d\mu$.
  In general this inequality is strict (Example \ref{exa:cantorian}).

  However we prove that $\|W\|_\mu^2 = \int_\calX \thet \,d\mu$ provided two additional conditions {\bf C1} and {\bf C2} hold (Proposition \ref{prop:dim=int}).

\item In Section \ref{sec:conv} we compute the $\mu$-norm of a convolution operator on $\calH = L^2(\mT)$, where $\mT=\mR/(2\pi\mZ)$ is a circle. For any bounded sequence $\{\lambda_k\}_{k\in\mZ}$ we consider the distribution
    $\lambda(x) = \sum \lambda_k e^{ikx}$. Then the operator
    $f\mapsto\Conv_\lambda f = \lambda * f$ is bounded:
    $\|\Conv_\lambda\| = \sup_{k\in\mZ} |\lambda_k|$. We prove (Proposition \ref{prop:piLambda}) that
    $$
      \|\Conv_\lambda\|_\mu^2 = \rho(\lambda), \qquad
      \rho(\lambda) = \limsup_{\# I\to\infty} \frac1{\# I} \sum_{k\in I} |\lambda_k|^2,
    $$
    where $I\subset\mZ$ are integer intervals.

\item In Section \ref{sec:dt} we compute $\mu$-norm for a wide class of bounded operators on $L^2(\mT)$, the so-called operators of diagonal type. Consider the operator $W = (W_{j,k})_{j,k\in\mZ}$
$$
  f = \sum f_k e^{ikx} \mapsto Wf = \sum_{j,k\in\mZ} W_{j,k} f_k e^{ijx}.
$$
We say that $W$ is of diagonal type ($W\in\DT$) if
\begin{equation}
\label{Wjk<c}
  |W_{j,k}| \le c_{j-k} \quad
  \mbox{for some sequence $\{c_s\}_{s\in\mZ}$ with $\sum c_s = \bc < \infty$}.
\end{equation}
We put $\|W\|_\DT = \inf\bc$ where infimum is taken over all sequences $c_s$, satisfying (\ref{Wjk<c}).

As simple examples we have the following operators of diagonal type.

(a) Convolution operators from Section \ref{sec:conv}.

(b) Operators of multiplication by functions with absolutely converging Fourier series.

(c) Linear combinations and products of operators of diagonal type.

(d) If $W\in\DT$ then the conjugated operator $W^*$ also lies in $\DT$.

We prove (Corollary \ref{cor:C*}) that the space $(\DT,\|\cdot\|_\DT)$ is a
  $C^*$-algebra.

\item We have the following inequality between norms ((\ref{|1|}) and
  Lemma \ref{lem:normdiagtype}):
$$
  \|W\|_\mu \le \|W\| \le \|W\|_\DT .
$$

\item We associate with $W\in\DT$ and any point $a\in\mT$ the distribution $L_a$,
$$
  L_a = \sum_{j\in\mZ} w_j(a) e^{ijx}, \qquad
  w_j(a) = \sum_{k\in\mZ} W_{j,k} e^{i(j-k)a}.
$$
We prove that

(1) the function $a\mapsto\rho(L_a)$ is continuous (Lemma \ref{lem:rhoLa}),

(2) $\displaystyle\|W\|_\mu^2 = \frac1{2\pi} \int_\mT \rho(L_a) \, da$ (Proposition \ref{prop:diagonal}).

\item For any operator $W\in\DT$ we introduce ``average trace'' of $W^* W$ by
$$
    \bT(W)
  = \limsup_{\# I\to\infty} \frac1{\# I} \sum_{j\in\mZ,\, l\in I} |W_{l,j}|^2,
$$
  where $I\subset\mZ$ are intervals. Then (Proposition \ref{prop:dimDT>})
  $\|W\|_\mu^2 \ge \bT(W)$.

We also prove (Proposition \ref{prop:T(UWU)}) that if $U$ is a unitary operator of diagonal type then $\bT(W) = \bT(U^{-1}WU)$.

\end{itemize}

\section{General properties of the $\mu$-norm}

\subsection{$\mu$-norm of $\pi_X$}

\begin{lem}
\label{lem:mupi}
For any $X\in\calB$
\begin{equation}
\label{pi=mu}
  \|\pi_X\|_\mu^2 = \mu(X).
\end{equation}
\end{lem}

{\it Proof}. For any partition $\chi = \{Y_1,\ldots,Y_J\}$
$$
    \calM_\chi(\pi_X)
  = \sum_{j=1}^J \mu(Y_j) \|\pi_X\pi_{Y_j}\|^2
  = \sum_{\mu(Y_j\cap X) > 0} \mu(Y_j)
 \ge \mu(X) .
$$
This implies $\calM_\chi(\pi_X) \ge \mu(X)$.

On the other hand, for $\chi = \{X,Y\}$, $Y=\calX\setminus X$ we have:
$$
    \calM_\chi(\pi_X)
  = \mu(X) \|\pi_X\|^2 + \mu(Y) \|\pi_X\pi_Y\|^2 = \mu(X).
$$
Hence $\inf_\chi \calM_\chi(\pi_X) = \mu(X)$.  \qed

\subsection{$\mu$-norm and subpartitions}

\begin{lem}
\label{lem:muM}
Let $W$ be a bounded operator and $\chi = \{Y_1,\ldots,Y_J\}$ a partition of the set $Y\in\calB$. Then
\begin{equation}
\label{subpartition}
  \mu(Y) \| W\pi_Y \|^2 \ge \sum_{j=1}^J \mu(Y_j) \| W\pi_{Y_j} \|^2.
\end{equation}
\end{lem}

{\it Proof}. The following two simple remarks
$$
  \|W\pi_{Y_j}\| \le \|W\pi_Y\| \quad
  \mbox{and}\quad
  \sum_{j=1}^J \mu(Y_j) = \mu(Y)
$$
imply (\ref{subpartition}).
 \qed

\begin{cor}
\label{cor:subpartition}
Lemma \ref{lem:muM} implies that $\calM_{\chi'}(W) \le \calM_{\chi}(W)$
if $\chi'$ is a subpartition of $\chi$.
\end{cor}

\subsection{Triangle inequality}
\label{sec:triangle}

\begin{lem}
\label{lem:tri}
For any two bounded operators $W_1$ and $W_2$
\begin{equation}
\label{trian0}
     \|W_1+W_2\|_\mu
 \le \|W_1\|_\mu + \|W_2\|_\mu .
\end{equation}
\end{lem}

{\it Proof}. First we show that for any partition $\chi=\{X_1,\ldots,X_J\}$
\begin{equation}
\label{trian00}
     \sqrt{\calM_\chi(W_1+W_2)}
 \le \sqrt{\calM_\chi(W_1)} +   \sqrt{\calM_\chi(W_2)} .
\end{equation}
We put
$$
  a_j = \|W_1\pi_{X_j}\|, \quad
  b_j = \|W_2\pi_{X_j}\|, \quad
  c_j = \|(W_1+W_2)\pi_{X_j}\|, \quad
  \mu_j = \mu(X_j).
$$
Then
\begin{equation}
\label{abc}
   0\le c_j\le a_j+b_j
\end{equation}
and inequality (\ref{trian00}) is equivalent to
$$
     \sqrt{\sum \mu_j c_j^2}
 \le \sqrt{\sum \mu_j a_j^2} +  \sqrt{\sum \mu_j b_j^2} .
$$
By (\ref{abc}) to prove this inequality, it is sufficient to check that
$$
     \sqrt{\sum \mu_j (a_j+b_j)^2}
 \le \sqrt{\sum \mu_j a_j^2} +  \sqrt{\sum \mu_j b_j^2} .
$$
This inequality is equivalent to
$$
     \Big(\sum \mu_j a_j b_j\Big)^2
 \le \sum \mu_j a_j^2 \cdot \sum \mu_j b_j^2
$$
which follows from Cauchy-Bunyakovsky-Schwartz.

To derive (\ref{trian0}) from (\ref{trian00}), we take a partition $\chi$ such that
$$
  \|W_1\|_\mu > \sqrt{\calM_\chi(W_1)} - \eps, \quad
  \mbox{and}\quad
  \|W_2\|_\mu > \sqrt{\calM_\chi(W_2)} - \eps .
$$
Then by (\ref{trian00})
\begin{eqnarray*}
       \|W_1+W_2\|_\mu
 &\le& \sqrt{\calM_\chi(W_1+W_2)}
  \le  \sqrt{\calM_\chi(W_1)} + \sqrt{\calM_\chi(W_2)} \\
 &\le& \|W_1\|_\mu + \|W_2\|_\mu + 2\eps.
\end{eqnarray*}
Since $\eps>0$ is arbitrary, (\ref{trian0}) follows. \qed

An equivalent form of inequality (\ref{trian0}) is
$$
      \big| \|W_2\|_\mu - \|W_1\|_\mu \big|
  \le \|W_2 - W_1\|_\mu .
$$

\begin{cor}
\label{cor:continuous}
By (\ref{WW}) we have:
$$
      \big| \|W_2\|_\mu - \|W_1\|_\mu \big|
  \le \|W_2 - W_1\| .
$$
Therefore, the $\mu$-norm is a continuous function of $W$ in the $L^2$ operator topology.
\end{cor}

\subsection{$U_{auto}(\calH)$-invariance}

Let $F:\calX\hookleftarrow$ be a measure preserving map: for any $X\in\calB$
$$
  F^{-1}(X) \in\calB, \quad\mbox{and}\quad
  \mu(X) = \mu(F^{-1}(X)).
$$
Such maps are called endomorphisms of the measure space $(\calX,\calB,\mu)$. If an endomorphism $F$ is invertible and the inverse map is measurable, then $F$ is said to be an automorphism.

Any endomorphism $F:\calX\hookleftarrow$ generates an operator
$U = U_F$ on $\calH$, where for any $f\in\calH$
\begin{equation}
\label{Uf=fF}
  U_F f = f\circ F.
\end{equation}
The operator (\ref{Uf=fF}) preserves the scalar product and satisfies the following identities:
\begin{eqnarray}
\label{id1}
     U_F(fg)
 &=& U_F f \cdot U_F g \qquad
     \mbox{for any }  f,g\in\calH,  \\
\label{id2}
     U_F \pi_\xi
 &=& \pi_{F^{-1}(\xi)} U_F \qquad
     \mbox{ for any }  \xi\in\calB.
\end{eqnarray}

Operators $U_F$ corresponding to automorphisms are invertible and therefore, unitary. They form a subgroup $U_{\auto}(\calH)$
in the group $U(\calH)$ of unitary operators on $\calH$.

Equation (\ref{id2}) implies the following

\begin{cor}
\label{cor:WU<W}
For any endomorphism $F$, any bounded $W$ and any $\xi\in\calB$
$$
  \|WU_F \pi_\xi\| \le \|W\pi_{F^{-1}(\xi)}\|.
$$

Indeed, for any $f\in L^2(\calX)$ (\ref{id2}) implies
$$
    \|WU_F \pi_\xi f\|
  = \|W \pi_{F^{-1}(\xi)} U_F f\|, \quad
    \| f \| = \| U_F f \|.
$$
\end{cor}

\begin{lem}
\label{lem:UWU}
For any bounded $W$ and any automorphism $F$
\begin{equation}
\label{invar}
  \|W U_F\|_\mu = \|W\|_\mu .
\end{equation}
\end{lem}

\begin{cor}
\label{cor:conj}
If $F$ is an automorphism then by (\ref{invar}) and (\ref{UWWU})
$\|U_F^{-1} W U_F\|_\mu = \|W\|_\mu$. Hence the $\mu$-norm is invariant w.r.t. conjugations by elements of $U_{\auto}(\calH)$.
\end{cor}

{\it Proof of Lemma \ref{lem:UWU}}. Let $\chi = \{Y_1,\ldots,Y_J\}$ be such that
$\|W\|_\mu^2 \ge \calM_\chi(W) - \eps$. Putting $\kappa = \{F(Y_1),\ldots,F(Y_J)\}$, we have:
$$
      \|WU_F\|_\mu^2
  \le \calM_\kappa(WU_F)
   =  \sum_{j=1}^J \mu(F(Y_j))  \|WU_F \pi_{F(Y_j)}\|^2 .
$$
By Corollary \ref{cor:WU<W} and the equations $F^{-1}(F(Y_j)) = Y_j$ this sum does not exceed
$$
     \sum_{j=1}^J \mu(Y_j) \|W\pi_{Y_j}\|^2
  =  \calM_\chi(W)
 \le \|W\|_\mu^2 + \eps.
$$
Since $\eps>0$ is arbitrary, $\|WU_F\|_\mu \le \|W\|_\mu$. Since $F$ is invertible, we have analogously $\|W\|_\mu \le \|WU_F\|_\mu$. \qed

\begin{cor}
\label{cor:piUpiU}
Equation (\ref{id2}) and Lemma \ref{lem:mupi} imply that for any collection of sets $Y_0,\ldots,Y_K\in\calB$ and any automorphism $F$
$$
     \| \pi_{Y_K} U_F \pi_{Y_{K-1}} U_F \ldots U_F \pi_{Y_1} U_F \pi_{Y_0} \|_\mu^2
  =  \mu\big( Y_K \cap \ldots \cap F^{-K+1}(Y_1) \cap F^{-K}(Y_0) \big) .
$$
\end{cor}

Indeed, by (\ref{id2}), (\ref{invar})
\begin{eqnarray*}
     \| \pi_{Y_K} U_F \pi_{Y_{K-1}} U_F \ldots \pi_{Y_1} U_F \pi_{Y_0} \|_\mu^2
 &=& \| \pi_{Y_K} \pi_{F^{-1}(Y_{K-1})} U_F^2 \pi_{Y_{K-2}}
              \ldots \pi_{Y_1} U_F \pi_{Y_0} \|_\mu^2 \\
 &=& \ldots
\; =\; \| \pi_{Y_K} \pi_{F^{-1}(Y_{K-1})} \ldots \pi_{F^{-K}(Y_0)} U_F^K \|_\mu^2 \\
 &=&  \| \pi_{Y_K \cap \ldots \cap F^{-K}(Y_0)} \|_\mu^2 .
\end{eqnarray*}
By (\ref{pi=mu}) the latter quantity equals
$\mu\big(Y_K \cap \ldots \cap F^{-k}(Y_0)\big)$.
\qed

\section{Examples}

\subsection{$\mu$-norm of a compact operator}

The measure $\mu$ is said to have no atoms if for any $X\in\calB$, $\mu(X)>0$ there exists a set $Y\in\calB$ such that $0<\mu(Y)<\mu(X)$.

\begin{lem}
\label{lem:finiterank}
Suppose $\mu$ has no atoms and $W$ is a finite rank operator. Then
$\|W\|_\mu = 0$.
\end{lem}

{\it Proof}. Any finite rank operator in $\calH$ has the form, (see, for example, \cite{Maurin})
$$
  f\mapsto Wf, \qquad
  Wf = \sum_{l=1}^L  \lambda_l \langle f,v_l\rangle u_l
$$
for two orthonormal systems $\{v_l\}$ and $\{u_l\}$. Given $\eps>0$ we can choose a partition
$\chi = \{Y_1,\ldots,Y_J\}$ with elements so small that
$$
  \|\pi_{Y_j} v_l\| \le \eps \quad
  \mbox{ for any $1\le j,l\le J$}.
$$
Then for any $f\in\calH$, $\|f\|=1$ and any $j\in \{1,\ldots,J\}$
$$
     \|W\pi_{Y_j} f\|^2
  =  \Big\| \sum_{l=1}^L  \lambda_l \langle\pi_{Y_j} f,v_l\rangle u_l \Big\|^2
 \le \Lambda\eps^2, \qquad
     \Lambda = \sum_1^L \lambda^2_l.
$$
Therefore $\calM_\chi(W) \le \sum \mu(Y_j) \Lambda\eps^2 = \Lambda \eps^2$. \qed

\begin{cor}
\label{cor:mu=0}
Suppose $\mu$ has no atoms and $W$ is a compact operator. Then $\|W\|_\mu =  0$.
\end{cor}

Indeed, any compact operator in $\calH$ may be approximated by finite rank operators. It remains to use Corollary \ref{cor:continuous}. \qed

\begin{cor}
\label{cor:B+K}
Suppose $\mu$ has no atoms. Then for any bounded $W$ and compact $W_0$ we have:
$\|W+W_0\|_\mu = \|W\|_\mu$.
\end{cor}

This statement is a combination of (\ref{trian0}) and Corollary \ref{cor:mu=0}. \qed

\subsection{Multiplication operator}

For any $g\in L^\infty(\calX)$ consider the operator $\widehat g$ of multiplication by $g$:
\begin{equation}
\label{W_g}
  L^2(\calX) \ni f \mapsto \widehat g f = gf.
\end{equation}

\begin{lem}
\label{lem:W_g}
$\|\widehat g\|_\mu = \|g\|$.
\end{lem}

{\it Proof}. First, suppose that $g$ is simple i.e., for some partition
$\kappa = \{X_1,\ldots,X_K\}$
$$
  g = \sum_{k=1}^K  c_k \one_{X_k}, \qquad
  c_k\in\mC.
$$
Then $\|\widehat g\pi_{Y}\|^2 = \max |c_k|^2$, where maximum is taken for $k$ such that
$\mu(Y\cap X_k) > 0$.

By Corollary \ref{cor:subpartition} it is sufficient to consider partitions
$\chi = \{Y_1,\ldots,Y_J\}$ which are subpartitions of $\kappa$. Then
$$
    \calM_\chi(\widehat g)
  = \sum_{j=1}^J \mu(Y_j) \max_{\mu(Y_j\cap X_k) > 0} |c_k|^2
  = \sum_{k=1}^K \mu(X_k) |c_k|^2
  = \|g\|^2.
$$

If $g$ is not simple then neglecting a set of a small measure, we can approximate $g$ in the $L^\infty$-norm by simple functions. \qed

\subsection{Finite-dimensional case}
\label{subsec:finite}

Let $\calX = \{1,\ldots,J\}$ be a finite set. Then $\calB$ consists of all subsets $X\subset\calX$. We assume that $\mu(X) = \#X / J$ i.e., the measure $\mu$ is uniformly distributed between points of $\calX$. Then $\calH = \mC^J$ with the product $\langle\,,\rangle = \frac1J (\,,)$, where $(\,,)$ is the standard Hermitian product.

Linear operators on $\calH$ are identified with $J\times J$-matrices. For any $X\subset\calX$
$$
  \pi_X = \diag(\delta_{1X},\ldots,\delta_{JX}), \qquad
    \delta_{jX}
  = \left\{\begin{array}{cc}
             1 \quad\mbox{if } j\in X, \\
             0 \quad\mbox{if } j\not\in X.
           \end{array}
    \right.
$$

\begin{lem}
\label{lem:finite}
Suppose $\#\calX=J<\infty$. Let $W = (W_{jk})$ be an operator on $\calH$. Then\footnote
{The quantity in the right-hand side of (\ref{finest+}) equals the trace of $W^* W$ divided by $J$. In particular, if $W$ is unitary then $\|W\|_\mu = 1$. Below such and analogous quantities will be said to be the average trace of $W^* W$.}
\begin{eqnarray}
\label{finest+}
      \|W\|_\mu^2
  &=& \frac1J \sum_{j,k=1}^J |W_{kj}|^2 , \\
\label{finest++}
      \|W\|_\mu
  &=& \|WU\|_\mu \quad
      \mbox{for any unitary operator $U$ on $\calH$}.
\end{eqnarray}
\end{lem}

{\it Proof}.  By Lemma \ref{lem:muM} to compute $\|W\|_\mu$, we may use the finest partition
\begin{equation}
\label{finest}
   \chi = \big\{ \{1\},\ldots,\{J\} \big\} .
\end{equation}
Then for any $j$ the unit vector $e_j = \pi_{\{ j\} } e_j$ has the coordinates
$e_{jk} = \sqrt{J} \delta_{jk}$. Therefore
$$
    \|W\|_\mu^2
  = \frac1J \sum_{j=1}^J \frac{\|W \pi_{\{j\}}\|^2}{J}
  = \frac1{J^2} \sum_{j=1}^J \|W e_j\|^2
  = \frac1J \sum_{j,k=1}^J |W_{kj}|^2 .
$$

To prove (\ref{finest++}) we may use (\ref{finest+}). Then
$$
    J\|WU\|_\mu^2
  = \sum_{j,k,l,s} W_{kl} U_{lj} \overline W_{ks} \overline U_{sj}
  = \sum_{j,k,l,s} W_{kl}\overline W_{ks} U_{lj} U^*_{js}
  = \sum_{k,l,s} W_{kl}\overline W_{ks} \delta_{ls}
  = \sum_{k,l} |W_{kl}|^2.
$$
The latter quantity equals $J\|W\|_\mu^2$.
\qed

\subsection{Dimension of a subspace}

Let $H\subset\calH$ be a closed subspace. We define its $\mu$-dimension by
$$
  \dim_\mu(H) := \|\pi_H\|_\mu^2 , \quad
  \mbox{where $\pi_H:\calH\to H$ is the orthogonal projector}.
$$

\begin{lem}
\label{lem:subspace}
1. If $\mu$ has no atoms and $H$ is finite-dimensional then $\dim_\mu(H)=0$.

2. If $H = \pi_X(\calH)$ for some $X\in\calB$ then $\dim_\mu(H) = \mu(X)$.

3. Suppose $\calX$ is finite, $\#\calX = J$, and for any $X\subset\calX$
$\mu(X) = \# X / J$. Then for any subspace $H\subset\calH=\mC^J$ its $\mu$-dimension equals its relative dimension: $\dim_\mu(H) = \dim(H) / J$, where $\dim (\cdot)$ is the conventional dimension of a vector space.
\end{lem}

{\it Proof}. The first statement follows from Lemma \ref{lem:finiterank}. The second statement follows from Lemma \ref{lem:mupi}. The third statement follows from Lemma \ref{lem:finite}. \qed

Suppose a cyclic group $\mZ_q$ acts on $(\calX,\calB,\mu)$ by automorphisms i.e., for any $s\in\mZ_q$ there is an automorphism $F_s$ such that the following identities hold:
\begin{equation}
\label{s'+s''}
  F_{s'} \circ F_{s''} = F_{s'+s''} \quad\mbox{and}\quad
  F_{-s} = F^{-1}_s.
\end{equation}
We assume that the action is almost free. This means that there is a set $X\in\calB$ such that the sets $X_s = F_s(X)$ do not intersect pairwise and
$\mu(\calX\setminus \cup X_s) = 0$. Let $U_s = U_{F_s}$ be the unitary operators, determined by (\ref{Uf=fF}). Then by (\ref{s'+s''})
$$
  U_{s'} U_{s''} = U_{s'+s''} \quad\mbox{and}\quad
  U_{-s} = U^{-1}_s.
$$

Putting $r=e^{2\pi i/q}$, we consider the spaces
$$
  H_n = \{f\in\calH : U_s f = r^{ns} f\;\mbox{ for any } s\in\mZ_q\}.
$$
Functions $f\in H_0$ are invariant w.r.t. the $\mZ_q$-action.

\begin{lem}
\label{lem:cyclic}
The spaces $H_n$ are pairwise orthogonal, $\calH = \oplus H_n$, and $\dim_\mu (H_n) = 1/q$.
\end{lem}

{\it Proof}. For any functions $f_j\in H_j$ and $f_k\in H_k$ and any $s\in\mZ_q$ we have:
$$
    \langle f_j,f_k\rangle
  = \langle U_s f_j, U_s f_k\rangle
  = r^{(j-k)s} \langle f_j,f_k\rangle.
$$
This proves the first statement of the lemma. The second statement follows from the equations ($f\in\calH$ is arbitrary)
\begin{equation}
\label{pi_Hn}
  f = \sum_{n\in\mZ_q} \pi_{H_n} f, \qquad
  \pi_{H_n} f = \frac1q \sum_{k\in\mZ_q} r^{-nk} U_k f \in H_n.
\end{equation}

To compute the $\mu$-dimension of $H_n$, consider a partition $\chi = \{Y_1,\ldots,Y_J\}$. We may assume that each element of the partition has nonempty intersection only with one set $X_s$. Then for any function $f = \pi_{Y_j} f$ the supports of the terms $r^{-nk} U_n f$ do not intersect pairwise. Therefore
$$
  \|\pi_{H_n} f\|^2 = \frac1{q^2} \sum_{k\in\mZ_q} \|r^{-nk} U_n f\|^2 = \frac{\|f\|^2}q.
$$
Hence $\calM_\chi(\pi_{H_n}) = \|f\|^2 / q$. This finishes the proof. \qed

\section{Additivity}
\label{sec:add0}

\begin{lem}
\label{lem:=}
For any partition $\kappa=\{X_1,\ldots,X_K\}$
\begin{equation}
\label{muW=}
     \|W\|_\mu^2
  =  \sum_{k=1}^K \|W\pi_{X_k}\|_\mu^2 .
\end{equation}
\end{lem}

{\it Proof}. Equation (\ref{muW=}) follows from the simple observation that for any subpartition $\chi = \{Y_1,\ldots,Y_J\}$ of the partition $\kappa$
$$
      \sum_{k=1}^K \calM_\chi(W\pi_{X_k})
   =  \calM_{\chi\vee\kappa}(W).
$$
\qed

It is natural to call the simple property of the $\mu$-norm presented in Lemma \ref{lem:=} the right additivity. The question about the left additivity turns out to be nontrivial. The following lemma presents a partial result in this direction.

\begin{lem}
\label{lem:<=}
For any partition $\kappa=\{X_1,\ldots,X_K\}$
\begin{equation}
\label{muW<=}
     \|W\|_\mu^2
 \le \sum_{k=1}^K \|\pi_{X_k} W\|_\mu^2 .
\end{equation}
\end{lem}

{\it Proof}. For any $f\in\calH$ and $Y\in\calB$
\begin{equation}
\label{orthogonal}
  \|W\pi_Y f\|^2 = \sum_{k=1}^K \|\pi_{X_k} W\pi_Y f\|^2,
\end{equation}
Therefore for any partition $\chi$
\begin{equation}
\label{M=M}
   \calM_\chi(W) \le \sum_{k=1}^K \calM_\chi(\pi_{X_k} W).
\end{equation}

\begin{rem}
Unfortunately, we do not have in (\ref{M=M}) an equation instead of the inequality. Indeed, for different terms in (\ref{orthogonal}) $\sup_{\|f\|=1}$ may be taken (almost taken) on different functions $f$.
\end{rem}

Let $\chi_1,\ldots,\chi_K$ be partitions such that
$$
  \calM_{\chi_k}(\pi_{X_k} W) \le \|\pi_{X_k} W\|_\mu^2 + \eps, \qquad
  k = 1,\ldots,K.
$$
We put $\chi = \vee_{k=1}^K \chi_k$. Then by (\ref{M=M}) and Corollary \ref{cor:subpartition}
$$
      \|W\|_\mu^2
  \le \calM_\chi(W)
  \le \sum_{k=1}^K \calM_\chi(\pi_{X_k} W)
  \le \sum_{k=1}^K \calM_{\chi_k}(\pi_{X_k} W)
  \le \sum_{k=1}^K \|\pi_{X_k} W\|_\mu^2 + K\eps.
$$
Since $\eps>0$ is arbitrary, we obtain (\ref{muW<=}). \qed

Let $\widehat g$ be the operator of multiplication by the function $g\in L^\infty(\calX)$.

\begin{lem}
\label{lem:Wr-additivity}
For any collection of functions $g_1,\ldots,g_K,g\in L^\infty(\calX)$ such that
\begin{equation}
\label{sumg=1}
  \sum |g_k|^2 = |g|^2
\end{equation}
and any bounded operator $W$
\begin{equation}
\label{Wg-add}
    \sum_{k=1}^K \|W\widehat g_k\|_\mu^2
  = \|W\widehat g\|_\mu^2 .
\end{equation}
\end{lem}

{\it Proof}. First, suppose that the functions $g_k$ are simple:
\begin{equation}
\label{g_simple}
   g_k = \sum_{j=1}^{J_k} \widetilde g_{kj} \one_{X_{kj}}, \qquad
   k = 1,\ldots,K,
\end{equation}
where $\widetilde g_{kj}\in\mC$ and $\chi_k = \{X_{k1},\ldots,X_{kJ_k}\}$ are partitions. Note that instead of the partitions $\chi_k$ we can use their common refinement $\chi = \chi_1\vee\ldots\vee\chi_K = \{X_1,\ldots,X_J\}$. Hence, we can replace (\ref{g_simple}) by
$$
   g_k = \sum_{j=1}^J g_{kj} \one_{X_j}, \qquad
   g_{kj}\in\mC.
$$
The function $g$ in this case is also simple: $g = \sum g_j\one_{X_j}$.
Equation (\ref{sumg=1}) implies
\begin{equation}
\label{g+g=g}
  \sum_{k=1}^K |g_{kj}|^2 = |g_j|^2 \quad
  \mbox{for any } j=1,\ldots,J.
\end{equation}

Equation (\ref{muW=}) implies
\begin{eqnarray*}
      \|W\widehat g\|_\mu^2
  &=& \sum_j \|W\widehat g \pi_{X_j}\|_\mu^2
   =  \sum_j \|W g_j\pi_{X_j}\|_\mu^2
   =  \sum_j |g_j|^2 \|W\pi_{X_j}\|_\mu^2 , \\
      \|W\widehat g_k\|_\mu^2
  &=& \sum_j \|W\widehat g_k \pi_{X_j}\|_\mu^2
   =  \sum_j \|W g_{kj}\pi_{X_j}\|_\mu^2
   =  \sum_j  |g_{kj}|^2 \|W\pi_{X_j}\|_\mu^2.
\end{eqnarray*}
Now (\ref{Wg-add}) follows from (\ref{g+g=g}).

If the functions $g_k$ are not simple, we approximate them in $L^\infty(\calX)$-norm by simple functions and use Corollary \ref{cor:continuous}. \qed

\section{An additional structure}
\label{sec:add}

Let the function $\dist(\cdot,\cdot):\calX\times\calX\to\mR_+$ determine on $\calX$ the structure of a compact metric space such that the balls
$$
  B_r(x) = \{y\in\calX : \dist(y,x) < r\}
$$
are measurable. Moreover, we assume that $\calB$ is the corresponding Borel $\sigma$-algebra while $\mu$ is a Borel measure.

\begin{prop}
\label{prop:theta}
For any $x\in\calX$ there exists the limit
\begin{equation}
\label{theta}
    \thet(x)
 =  \lim_{\eps\searrow 0} \thet_\eps(x), \qquad
    \thet_\eps(x)
 =  \|W\pi_{B_\eps(x)}\|^2,
\end{equation}
the function $\thet:\calX\to [0,\|W\|^2]$ is measurable and
\begin{equation}
\label{dim<inttheta}
  \|W\|_\mu^2 \le \int_\calX \thet\, d\mu.
\end{equation}
\end{prop}

{\it Proof of Proposition \ref{prop:theta}}. For any $x\in\calX$ the function $\eps\mapsto\thet_\eps(x)$, $\eps>0$ is non-decreasing and non-negative. This implies the existence of limits (\ref{theta}).

Consider a sequence $\{\eps_j\}_{j\in\mN}$, $\eps_j\to 0$ as $j\to\infty$. For any $j\in\mN$ let $\calN_j = \{y_1,\ldots,y_K\} \subset\calX$ be a finite $\eps_j$-net\footnote
{The number $K\in\mN$ as well as the points $y_k$ depend on $j$, but for brevity we do not indicate this in the notation.}:
$$
  \mbox{for any $x\in\calX$ there exists $y_k\in\calN_j$ such that $\dist(y_k,x) < \eps_j$}.
$$
We define
$$
      y(x,j)
  \in \{y_k\in\calN_j : \dist(y_k,x) = \min_{y\in\calN_j} \dist(y,x)\}
$$
as the element with minimal possible index. Then $\dist(y(x,j),x) < \eps_j$.

For any $x\in\calX$ and $j\in\mN$
\begin{equation}
\label{BBB}
  B_{\eps_j}(x) \subset B_{2\eps_j}(y(x,j)) \subset B_{3\eps_j}(x).
\end{equation}
We put $\widetilde\thet_j(x) = \|W\pi_{B_{2\eps_j}(y(x,j))}\|^2$. Then $\widetilde\thet_j$ takes a finite number of values and by construction of $y(x,j)$ the preimage of each value is measurable. Therefore $\widetilde\thet_j$ is measurable.

Inclusions (\ref{BBB}) imply
\begin{equation}
\label{ttt}
  \thet_{\eps_j}(x) \le \widetilde\thet_j(x) \le \thet_{3\eps_j}(x) .
\end{equation}
Inequalities (\ref{ttt}) imply
$$
    \lim_{j\to\infty} \widetilde\thet_j(x)
  = \lim_{\eps\searrow 0} \thet_\eps(x)
  = \thet(x) .
$$
Hence $\thet$ is measurable as a pointwise limit of measurable functions, \cite{Dudley}.
\smallskip

Now we turn to the proof of inequality (\ref{dim<inttheta}).
We may assume that $\eps_{j+1} \le \eps_j/3$. Then by (\ref{BBB})
$B_{2\eps_{j+1}}(y(x,j+1))\subset B_{2\eps_j}(y(x,j))$. Therefore $\widetilde\thet_{j+1}(x)\le\widetilde\thet_j(x)$ for any $x\in\calX$. The sequences $\{\widetilde\thet_j(x)\}_{j\in\mN}$ monotonically tend to $\thet(x)$ as $j\to\infty$. The sets
$$
  S_{j,\sigma} = \{x\in\calX : \widetilde\thet_j - \thet(x) > \sigma\}
$$
are measurable and satisfy for any $\sigma>0$
\begin{equation}
\label{S=}
  S_{j+1,\sigma} \subset S_{j,\sigma}, \quad
  \lim_{j\to\infty} \mu(S_{j,\sigma}) = 0.
\end{equation}

We put $\eps' = \eps_{j'} / 2$, where $\mu(S_{j',\sigma}) < \sigma_1$.
By (generalization of) Lusin's theorem \cite{Feldman} there exists $X_{\sigma_1}\in\calB$ such that
$\mu(X_{\sigma_1}) \le \sigma_1$, $\calX_{\sigma_1} = \calX\setminus X_{\sigma_1}$ is compact, and $\thet$ is continuous on $\calX_{\sigma_1}$.

The function $\thet|_{\calX_{\sigma_1}}$ is uniformly continuous. Hence there exists $\eps''>0$ such that for any $x,y\in\calX_{\sigma_1}$
\begin{equation}
\label{thet-thet}
  |x - y| < \eps'' \quad\mbox{implies}\quad
  |\thet(x) - \thet(y)| < \sigma_2.
\end{equation}

For $0<\eps<\min\{\eps',\eps''\}$ let $\calN_\eps = \{y_1,\ldots,y_K\} \subset\calX_{\sigma_1}$ be a finite $\eps$-net on $\calX_{\sigma_1}$. The balls $B_\eps(y_1),\ldots,B_\eps(y_K)$ cover $\calX_{\sigma_1}$. Consider the partition $\chi = \{S,X_1,\ldots,X_K\}$, where
$$
  S = X_{\sigma_1}\cup S_{j',\sigma}, \quad
  X_1 = B_\eps(y_1)\setminus S, \quad
  X_j = B_\eps(y_j)\setminus (S\cup B_\eps(y_1) \cup\ldots,\cup B_\eps(y_{j-1}))
$$
for $j\ge 2$. Then
$$
    \Theta(x)
  = \left\{ \begin{array}{ccc}
             \thet(y_j) & \mbox{if} & x\in X_j, \\
                 0      & \mbox{if} & x\in S
            \end{array} \right.
$$
is a simple function on $\calX$ which approximates $\thet$ on the set
$\calX\setminus S$ with precision $\sigma_2$:
$$
      \big| (\Theta - \thet)|_{\calX\setminus S} \big|
  \le \sigma_2.
$$
Moreover,
\begin{equation}
\label{intS}
      \Big| \int_\calX \Theta\, d\mu - \int_\calX \thet\, d\mu \Big|
  \le \sigma_2 + \int_S \thet\, d\mu
  \le \sigma_2 + 2\|W\|^2 \sigma_1.
\end{equation}

Consider the quantity
$\displaystyle
     \calM_\chi(W)
  =  \mu(S) \|W\pi_S\|^2
   + \sum_{k=1}^K \mu(X_k) \|W\pi_{X_k}\|^2$.
Then
\begin{equation}
\label{M-Theta}
     \calM_\chi(W) - \int_\calX \Theta\, d\mu
 \le 2\sigma_1 \|W\|^2 + \sum_{k=1}^K \mu(X_k) (\|W\pi_{X_k}\|^2 - \thet(y_k)).
\end{equation}

By (\ref{BBB}) $X_k\subset B_{2\eps_{j'}}(y(y_k,j'))$. Therefore by (\ref{thet-thet})
$$
  \|W\pi_{X_k}\|^2 \le \thet(y_k) + \sigma_2.
$$
This implies $\sum \mu(X_k) \|W\pi_{X_k}\|^2 \le \int_\calX \Theta\, d\mu + \sigma_2$. Combining this estimate with (\ref{intS}) and (\ref{M-Theta}), we obtain:
$$
      \calM_\chi(W)
  \le \int_\calX \thet\,d\mu + 2\sigma_2 + 2\|W\|^2 \sigma_1.
$$
Since $\sigma_1$ and $\sigma_2$ may be chosen arbitrary small, we obtain (\ref{dim<inttheta}). \qed

\begin{exa}
\label{exa:cantorian}
Inequality in (\ref{dim<inttheta}) may be strict. Indeed, let $\calX$ be the circle $\mT=\mR/(2\pi\mZ)$ with the Lebesgue measure $d\mu = dx/(2\pi)$. Let $X\subset\mT$ be the complement to a Cantorian set of positive measure. Then $X$ is open, dense, and $\mu(X)<1$. Take $W = \pi_X$. Then $\thet\equiv 1$. However
$$
  \|W\|_\mu^2 = \mu(X) < 1 = \frac1{2\pi} \int_\mT \thet(x)\, dx.
$$
\end{exa}

To replace in (\ref{dim<inttheta}) the inequality by the equation, we need additional assumptions.
\medskip

{\bf C1}. The function $\thet$ is continuous.
\medskip

{\bf C2}. There exists $c>0$ such that for any open $\calO\subset\calX$, $\diam(\calO)\le\eps$ and any $x\in\calO$ there exists $f = \pi_\calO f$ satisfying
$$
      \Big| \|Wf\|^2 - \thet(x) \|f\|^2 \Big|
  \le \gamma(\eps) \|f\|^2, \quad
      c < \big| f|_\calO \big| < c^{-1},
$$
where $\gamma(\eps)\to 0$ as $\eps\to 0$.

\begin{prop}
\label{prop:dim=int}
Suppose conditions {\bf C1}--{\bf C2} hold. Then
$$
  \|W\|_\mu^2 = \int_\calX \thet\,d\mu.
$$
\end{prop}

{\it Proof}. According to Proposition \ref{prop:theta} it is sufficient to prove the inequality
$$
  \|W\|_\mu^2 \ge \int_\calX \thet \, d\mu .
$$
By {\bf C1} $\thet$ is continuous on $\calX$. Therefore it is uniformly continuous. Given a small $\sigma>0$ we choose $\eps>0$ such that for any ball $B\subset\calX$ of radius $\eps$
\begin{equation}
\label{thet_continuous}
  x,y\in B \quad \mbox{implies}\quad
  |\thet(x) - \thet(y)| \le \sigma.
\end{equation}

Suppose
$$
  \|W\|_\mu^2 < \int_\calX \thet\,d\mu - \alpha \quad
  \mbox{for some $\alpha>0$}.
$$
Then for some partition $\chi=\{Y_1,\ldots,Y_K\}$ we have:
$$
  \calM_\chi(W) < \int_\calX \thet\,d\mu - \alpha.
$$
We are going to show that this inequality holds only for $\alpha=\alpha(\eps)\to 0$ as $\eps\to 0$.

By Lemma \ref{lem:muM} we may assume that $\diam(Y_k) < \eps$, $k=1,\ldots,K$. We take arbitrary $Y\in\chi$ such that
\begin{equation}
\label{thet-alpha}
     \mu(Y) \|W\pi_Y\|^2 < \int_Y \thet \, d\mu - \alpha.
\end{equation}
By \cite{Dudley} the measure $\mu$ is regular as any Borel probability measure on a compact metric space.
Hence for any $\gamma_1>0$ there exists a ball $B_\eps(x)$ and an open set $A$ such that
\begin{equation}
\label{muA}
  Y\subset A\subset B_\eps(x), \quad
  \mu(Y) \ge \mu(A) (1 - \gamma_1).
\end{equation}
Putting $D = A\setminus Y$, we have by (\ref{muA}):
\begin{equation}
\label{muD<muA}
  \mu(D) \le \gamma_1 \mu(A).
\end{equation}

Let $f = \pi_A f$ be a function which by {\bf C2} satisfies
\begin{equation}
\label{x'}
      \Big| \|Wf\|^2 - \thet(x) \|f\|^2 \Big|
  \le \gamma(\eps) \|f\|^2 , \quad
      c < \big|f|_A\big| < c^{-1} , \quad
      x' \in A.
\end{equation}
Since $f = \pi_A f = \pi_{A\cap Y} f + \pi_D f$, the triangle inequality implies
$$
       \|W\pi_Y f\|
  \ge  \|W f\| - \|W\pi_D f\|
$$
We estimate $\|Wf\|$ by (\ref{x'}) and $\|W\pi_D f\| \le \|W\|\,\|\pi_D f\|$ by Lemma \ref{lem:piDf}:
$$
  \|Wf\| \ge (\thet(x) - \gamma(\eps)) \|f\|^2 , \quad
  \|\pi_D f\| \le \frac{\sqrt{\mu(D)}}{c^2 \sqrt{\mu(A)}} \|f\| .
$$
Therefore for any $x'\in A\cap Y$
\begin{eqnarray*}
       \|W\pi_Y f\|^2
 &\ge& \bigg( \sqrt{\thet(x') - \gamma(\eps)} \|f\|
            - \|W\| \frac{\sqrt{\mu(D)}}{c^2\sqrt{\mu(A)}} \|f\| \bigg)^2 \\
 &\ge&  \bigg( \thet(x') - \gamma(\eps) - 2\frac{\sqrt{\thet(x)} \|W\|}{c^2}
                                              \frac{\sqrt{\mu(D)}}{\sqrt{\mu(A)}}
        \bigg) \|f\|^2.
\end{eqnarray*}

By using (\ref{thet_continuous}), (\ref{muD<muA}), and the estimate $\thet \le \|W\|^2$, we obtain
$$
      \frac{\|W\pi_Y f\|^2}{\|f\|^2}
  \ge \thet(x) - \sigma - \gamma(\eps) - 2 \sqrt{\gamma_1} \frac{\|W\|^2}{c^2} .
$$

This estimate implies that
$$
      \|W\pi_Y\|^2 - \frac1{\mu(Y)}\int_Y\thet\, d\mu
  \ge \frac{\|W\pi_Y f\|^2}{\|f\|^2} - \frac1{\mu(Y)}\int_Y\thet\, d\mu
  \ge - \alpha,
$$
where
\begin{eqnarray*}
     \alpha
 &=& \Big| \thet(x) - \frac1{\mu(Y)} \int_Y \thet\,d\mu \Big|
    + \gamma_1\thet(x) + \sigma + \gamma(\eps)
    + 2\sqrt{\gamma_1} \frac{\|W\|^2}{c^2} \\
 &=&  2\sigma + \gamma_1\thet(x) + \gamma(\eps)
    + 2\sqrt{\gamma_1} \frac{\|W\|^2}{c^2} .
\end{eqnarray*}
Hence $\alpha$ in (\ref{thet-alpha}) has to be arbitrarily small if $\eps$ is small enough. \qed

\section{Convolutions on $L^2(\mT)$}
\label{sec:conv}
\subsection{Preliminary constructions}

{\bf a}. Let $\mT = \mR/(2\pi\mZ)$ be a circle. The measure $\mu$ and the Hermitian product have the form
$$
  d\mu = \frac1{2\pi}\,dx, \quad
  \langle f,g\rangle = \frac1{2\pi} \int_\mT f\overline g \, dx.
$$

For any bounded sequence $\{\lambda_k\}_{k\in\mZ}$ consider the distribution
\begin{equation}
\label{lambda(x)}
  \lambda(x) = \sum_{k\in\mZ} \lambda_k e^{ikx}.
\end{equation}
In particular, if all $\lambda_k$ equal $1$ then $\lambda=2\pi\delta$, where $\delta=\delta(x)$ is the $\delta$-function on $\mT$.

{\bf b}. The operator of convolution
$$
  L^2(\mT) \ni f \mapsto \Conv_\lambda f = \lambda * f = \frac1{2\pi}\int_{\mT} \lambda(x-y) f(y)\, dy
$$
is bounded and
\begin{equation}
\label{Cnorm}
  \|\Conv_\lambda\| = c_\lambda, \qquad c_\lambda =  \sup_{k\in\mZ} |\lambda_k| .
\end{equation}

\begin{rem}
\label{rem:priority}
Below to reduce the number of brackets, we assume that the convolution has a higher priority in comparison with the arithmetic operations. Hence, $f g * h$ means $f\cdot (g * h)$ for any functions $f,g,h$ on $\mT$.
\end{rem}

If the sequence $\{\lambda_k\}_{k\in\mZ}$ lies in the space $l^2$ then $\lambda\in L^2(\mT)$ and
\begin{equation}
\label{|d|}
  \|\lambda\|^2 = \|\{\lambda_k\}_{k\in\mZ}\|^2_{l^2} .
\end{equation}

{\bf c}. For any $n\in\mZ$, any $a\in L^\infty(\mT)$ and any $g\in L^2(\mT)$
\begin{equation}
\label{d*e}
  \lambda * (e^{inx} g) = e^{inx} ( (e^{-inx}\lambda) * g), \quad
  \| ag \| \le \|a\|_\infty \|g\| ,
\end{equation}
where $\|\cdot\|_\infty$ is the norm in $L^\infty$.
\smallskip

{\bf d}. For any $\Lambda\subset\mZ$ consider the sequence $\{\lambda_k\}_{k\in\mZ}$, where $\lambda_k=1$ if $k\in\Lambda$ and $\lambda_k = 0$ if $k\not\in\Lambda$. Let $d_\Lambda$ denote the corresponding distribution (\ref{lambda(x)}). Then
$\Conv_{d_\Lambda}$ is the orthogonal projector to the subspace of $L^2(\mT)$ spanned by the vectors $e^{ikx}$, $k\in\Lambda$.

Given $a\in\mT$ consider the distribution $d_{a,\Lambda}$, defined by
$$
  d_{a,\Lambda}(x) = d_{\Lambda}(x-a)\quad  \mbox{for any } x\in\mT.
$$

If the sets $A,B\subset\mZ$ have empty intersection then $d_{a,A\cup B} = d_{a,A} + d_{a,B}$, and for any $a\in\mT$ and any $g\in L^2(\mT)$
\begin{equation}
\label{d+d}
  \|d_{a,A\cup B} * g\|^2 = \|d_{a,A} * g\|^2 + \|d_{a,B} * g\|^2 .
\end{equation}

For any $k\in\mZ$ let $\{k\}$ denote the corresponding one-point subset of $\mZ$. Then for any $f\in L^2(\mT)$
\begin{equation}
\label{Fcoeff}
  d_{a,\{k\}} * f = f_k d_{a,\{k\}} .
\end{equation}
where $f_k$ is the Fourier coefficient of $f$ with the number $k$.

{\bf e}. For any interval $J = \{n,n+1,\ldots,m\}\subset\mZ$ let $\# J = m-n+1$ be the number of points on $J$.

\begin{lem}
\label{lem:YfJ}
Let $Y=[a-\eps,a+\eps]$ and
$\displaystyle   f = \pi_Y f = \sum_{k\in\mZ} f_k e^{ikx} \in L^2(\mT)$.
Then for any integer interval $J$ and any $m\in\mZ$
\begin{eqnarray}
\label{f-ef}
       \|f - e^{im(x-a)} f\|
 &\le& |m| \eps \|f\| , \\
\label{f-fe}
       |f_m - e^{ila} f_{m+l}|
 &\le& \frac{\eps^{3/2}}{\sqrt\pi} |l| \|f\| , \\
\label{ff-f}
       \Big| \sum_{k\in\mZ} e^{-ima} f_k\overline f_{k+m} - \|f\|^2 \Big|
 &\le& |m| \eps \|f\|^2 .
\end{eqnarray}
\end{lem}

{\it Proof}. We put $g = f - e^{im(x-a)} f$. By using inequality (\ref{d*e}), we have:
$$
    \|g\|^2
  = \frac1{2\pi} \int_{a-\eps}^{a+\eps} |1 - e^{im(x-a)}|^2 |f|^2\, dx
 \le \|\pi_Y (1 - e^{im(x-a)})\|_\infty^2 \|f\|^2
 \le m^2 \eps^2 \|f\|^2.
$$
This implies (\ref{f-ef}).

Now we turn to inequality (\ref{f-fe}). We have:
\begin{eqnarray*}
      |f_m - f_{m+l} e^{ila}|
  &=& \bigg| \frac1{2\pi} \int_{a-\eps}^{a+\eps}
          \Big( f(x) e^{-ilx} - f(x) e^{-ilx-il(x-a)}
          \Big)\, dx \bigg|  \\
 &\le& \frac1{2\pi} \int_{a-\eps}^{a+\eps}
          \Big| f(x) (1 - e^{-il(x-a)}) \Big|\, dx \\
 &\le& \|f\|\, \big\|\pi_Y(1 - e^{-il(x-a)})\big\|
 \;\le\; \frac{\eps^{3/2}}{\sqrt\pi} |l|\, \|f\| .
\end{eqnarray*}

To prove (\ref{ff-f}), we consider the Hermitian product
$$
  \langle g,f\rangle = \|f\|^2 - \sum_{k\in\mZ} e^{-ima} f_k \overline f_{k+m} .
$$
By (\ref{f-ef}) $|\langle g,f\rangle| \le \|g\|\cdot\|f\| \le |m|\eps \|f\|^2$.
Hence, (\ref{ff-f}) follows. \qed

{\bf f}. For any distribution (\ref{lambda(x)}) we put
\begin{equation}
\label{rho}
  \rho(\lambda) = \limsup_{\# I\to\infty} \rho_I(\lambda), \qquad
  \rho_I(\lambda) = \frac{1}{\# I} \sum_{k\in I} |\lambda_k|^2 ,
\end{equation}
where $I$ denotes integer intervals.

\begin{prop}
\label{prop:piLambda}
$\|\Conv_\lambda\|_\mu^2 = \rho(\lambda)$.
\end{prop}

\begin{rem}
Corollary \ref{cor:mu=0} suggests the question: does the equation $\|W\|_\mu=0$ imply compactness of the operator $W$? Proposition \ref{prop:piLambda} gives an example of a noncompact operator $W$ with $\|W\|_\mu=0$. We may take $W=\Conv_{d_\Lambda}$, where $\Lambda=\{k^2 : k\in\mN\}$.
\end{rem}

We prove Proposition \ref{prop:piLambda} in Sections \ref{sec:comp} and \ref{sec:appl}

\subsection{Computation of $\|\Conv_\lambda\pi_{(-\eps,\eps)}\|$}
\label{sec:comp}

We put $Y = (a-\eps,a+\eps)\subset\mT$. It is possible to change the variable $x\mapsto x-a$ on $\mT$. This will change $\lambda(x)$ by $\lambda(x-a)$, but $\rho$ will remain the same. Hence we can assume that $a=0$ and
$Y = (-\eps,\eps)$.

\begin{lem}
\label{lem:conv}
For any $\gamma>0$ and any $f = \pi_Y f\in L^2(\mT)$
\begin{equation}
\label{d*f}
  \|\lambda * f\|^2 \le (\rho + \gamma)\|f\|^2  \quad
  \mbox{if $\eps$ is sufficiently small.}
\end{equation}
For any $\gamma>0$ there exists $f = \pi_Y f$, such that $|f|=\mbox{\rm const}$ on $Y$ and
\begin{equation}
\label{inverse}
    \|\lambda * f\|^2 \ge (\rho - \gamma)\|f\|^2 \quad
    \mbox{if $\eps$ is sufficiently small.}
\end{equation}
\end{lem}

\begin{cor}
\label{cor:conv}
Estimates (\ref{d*f})--(\ref{inverse}) imply that for any sufficiently small interval $Y$
\begin{equation}
\label{|conv|=o}
  \Big| \|\Conv_\lambda \pi_Y\| - \rho \Big| = o(\mu(Y)).
\end{equation}
\end{cor}

{\it Proof of Lemma \ref{lem:conv}}.
Let $M'(\sigma)$ be such that
\begin{equation}
\label{M'}
  \rho_I(\lambda) \le \rho + \sigma \quad
  \mbox{for any interval $I\subset\mZ$, \quad $\# I \ge M'(\sigma)$}.
\end{equation}

Let $I=-I\subset\mZ$ be an interval with $\# I = 2B + 1 > M'(\sigma)$. We put
\begin{equation}
\label{F=}
  F = \frac1{\# I} d_I f
    = \frac1{2B+1} \sum_{k\in\mZ,\, |l|\le B} f_{k+l} e^{ikx}.
\end{equation}
By Lemma \ref{lem:YfJ}
$$
      \|f - F\|
  \le \frac1{\# I} \sum_{l\in I} \|f - e^{ilx} f\|
  \le \frac1{2B+1} \sum_{|l|\le B} |l|\eps \|f\|.
$$
Hence
\begin{equation}
\label{f-F}
  \|f - F\| \le (B+1)\eps \|f\|, \quad
  \|\lambda * f - \lambda * F\| \le c_\lambda (B+1)\eps \|f\|
\end{equation}

We have:
\begin{eqnarray}
\nonumber
     \|\lambda * F\|^2
 &=& \frac1{(2B+1)^2} \sum_{k\in\mZ} |\lambda_k|^2 \sum_{|l|,|n|\le B} f_{k+l}\overline f_{k+n} \\
\nonumber
 &=& \frac1{(2B+1)^2} \sum_{k\in\mZ} |\lambda_{k-l}|^2 \sum_{|l|,|m+l|\le B} f_k\overline f_{k+m} \\
\label{|lam*F|}
 &=& \sum_{|m|\le 2B} \sum_{k\in\mZ} b_m \rho_{I_{k,m}} f_k\overline f_{k+m},
\end{eqnarray}
where
\begin{equation}
\label{Ikm}
    b_m
  = \frac{2B + 1 - |m|}{(2B+1)^2}, \quad
    I_{k,m}
  = \left\{ \begin{array}{cc}
              [k-B,k+B-m] \cap\mZ &  \mbox{ if } m\ge 0, \\ {}
              [k-B-m,k+B] \cap\mZ &  \mbox{ if } m < 0.
            \end{array}
    \right.
\end{equation}
Note that
\begin{equation}
\label{sumb}
  \sum_{|m|\le 2B} b_m = 1, \quad
  \# I_{k,m} = 2B + 1 - |m|, \quad
  \rho_{I_{k,m}} \le c_\lambda,
\end{equation}
and by (\ref{M'}) $\rho_{I_{k,m}} \le \rho+\sigma$ if $\# I_{k,m} \ge M' = M'(\sigma)$.

By using (\ref{ff-f}) with $a=0$, we have the estimate
\begin{eqnarray*}
       \|\lambda * F\|^2
 &\le& \sum_{2B+1-|m|<M',\,k\in\mZ} b_m c_\lambda f_k\overline f_{k+m}
      + \sum_{2B+1-|m|\ge M',\,k\in\mZ} b_m (\rho - \sigma) f_k\overline f_{k+m} \\
 &\le& \bigg( \frac{M'(1+M')}{(2B+1)^2} c_\lambda
             + \frac{(2B+1)^2 - M'(1+M')}{(2B+1)^2} (\rho - \sigma) \bigg)
              (1 - B\eps) \|f\|^2.
\end{eqnarray*}
This implies (\ref{d*f}) if $\sigma$, $M'/B$, and $B\eps$ are small.
\smallskip

Now we check estimate (\ref{inverse}).
We fix small positive $\sigma$ and $\eps$.
By Lemma \ref{lem:C<J<2C} we choose an interval $J\subset\mZ$ such that
\begin{eqnarray}
\label{J1}
    \# J \ge 20M'(\sigma)   \mbox{ is odd}, \quad
    \rho_J(\lambda)  >  \rho - \sigma, \\
\label{J2}
   2(\sigma_1\eps)^{-1} + 1 \le \# J \le 4(\sigma_1\eps)^{-1} + 3, \quad
   \sigma_1 = \sigma^{1/3}.
\end{eqnarray}

Consider the function
$$
  \ph = \pi_Y \ph = \frac{\one_Y}{\sqrt{2\eps}}
      = \frac1{\sqrt{2\eps}} \sum_{k\in\mZ} \frac{\sin k\eps}{\pi k} e^{ikx}, \quad
  \|\ph\| = 1.
$$
For any interval $I=-I\subset\mZ$ such that $\# I = 2b+1$ we have the estimate
\begin{equation}
\label{f0}
      \|\ph - d_I * \ph\|^2
  \le \frac1{2\eps} \sum_{|k|>b} \frac1{\pi^2 k^2}
   <  \frac1{b\eps} .
\end{equation}

By (\ref{J1}) $\# J$ is odd. Hence there exists $q\in\mZ$ such that $(J-q)=-(J-q)$. We define $f=e^{iqx}\ph$ and show that this $f$ satisfies (\ref{inverse}).

Obviously $f=\pi_Y f$ and $|f|=\mbox{const}$ on $Y$. We put
$$
  I = [-B,B]\cap\mZ, \quad
  I' = [-M',M']\cap\mZ, \quad
  M' = M'(\sigma), \quad
  B > M'.
$$
and consider $\widetilde J = -\widetilde J\subset\mZ$ such that $\widetilde J + 4I + I' = J-q$.
We will assume
\begin{equation}
\label{8B+2M}
  (\sigma_1\eps)^{-1} > 8B + 2M'.
\end{equation}
Then by (\ref{J2})
$$
    \frac1{\eps\sigma_1}
  < \#\widetilde J = \# J - 8B - 2M'
  < \frac4{\eps\sigma_1}
$$
and by (\ref{f0}) with $b=(\sigma_1\eps)^{-1}$
\begin{eqnarray}
\label{f-d*f}
       \|f - d_{\widetilde J + q} * f\|^2
 &\le& \sigma_1 \|f\|^2 ,  \\
\label{df-ddf}
       \|\lambda * f - \lambda * d_{\widetilde J + q} * f\|
 &\le& c_\lambda \|f - d_{\widetilde J + q} * f\|
  \le  c_\lambda \sqrt{\sigma_1} \|f\|.
\end{eqnarray}

We determine $F$ by (\ref{F=}) and put
\begin{equation}
\label{Phi}
    \Phi
  = \sum_{|m|\le 2B} \sum_{k\in\widetilde J + 2I + q} b_m \rho_{I_{k,m}} f_k\overline f_{m+k}.
\end{equation}
Since $\rho_I\le c_\lambda$, for any interval $I\subset\mZ$, we have by (\ref{|lam*F|}):
\begin{eqnarray*}
      \Big| \Phi - \|\lambda * F\|^2 \Big|
 &=&  \sum_{|m|\le 2B} \sum_{l\not\in\widetilde J + 2I} b_m \rho_{I_{l+q,m}} \ph_l\overline \ph_{m+l} \\
&\le&  \sum_{|m|\le 2B} \sum_{l\not\in\widetilde J + 2I} b_m c_\lambda |\ph_l\ph_{m+l}| \\
&\le&  \sum_{|m|\le 2B} \sum_{l\not\in\widetilde J} c_\lambda b_m |\ph_l|^2
\;\le\; 2c_\lambda \sigma_1 \|f\|^2.
\end{eqnarray*}

Note that the conditions $|m|\le 2B$, $k\in\widetilde J + 2I + q$ in (\ref{Phi}) imply $I_{k,m}\subset J$, where $I_{k,m}$ are defined in (\ref{Ikm}). Therefore by Lemma \ref{lem:rhoJ0}
$$
  \rho_{I_{k,m}} = \rho - s_m, \qquad
  s_m \le \frac{2\sigma \# J}{2B+1-|m|}.
$$
This implies that $\Phi$ is estimated as follows:
$$
       \Phi
  \ge  \sum_{|m|\le 2B} \sum_{k\in\widetilde J + 2I + q}
         \Big( b_m\rho - \frac{2\sigma \# J}{(2B+1)^2} \Big) f_k\overline f_{m+k}
$$
We perform summation in $k$ with the help of (\ref{ff-f}) and (\ref{f-d*f}):
\begin{eqnarray*}
      \sum_{k\in\widetilde J + 2I + q} f_k\overline f_{m+k}
 &=&  \bigg( \sum_{k\in\mZ} \;\; - \!\! \sum_{k\not\in\widetilde J + 2I + q} \bigg) f_k\overline f_{m+k} \\
&\ge&  (1 - |m| \eps) \|f\|^2 - \sum_{k\not\in\widetilde J + q} |f_k|^2
\;\ge\;  (1 - |m| \eps - \sigma_1) \|f\|^2 .
\end{eqnarray*}
Hence
\begin{eqnarray*}
       \Phi
 &\ge& \sum_{|m|\le 2B} (1 - 2\eps B - \sigma_1)
         \Big( b_m\rho - \frac{2\sigma \# J}{(2B+1)^2} \Big) \|f\|^2 \\
 &\ge& (1 - 2\eps B - \sigma_1) \Big( \rho - \frac{2\sigma \# J}{2B+1}  \Big) \|f\|^2 .
\end{eqnarray*}

This estimate implies that
$$
  \Phi \ge (\rho - \gamma) \|f\|^2,
$$
where $\gamma$ is small provided $\sigma_1$, $\sigma\# J / B$, and $\eps B$ are small. These quantities should be small, and conditions (\ref{J2}) and (\ref{8B+2M}) should hold. This can be made if we choose
$\sigma_1 = \sigma^{1/3}$, $\# J / 11 < B < \# J / 10$. \qed

\subsection{Application of Proposition \ref{prop:dim=int}}
\label{sec:appl}

Corollary \ref{cor:conv} implies that the function $\thet$, corresponding by (\ref{theta}) to the operator $W=\Conv_\lambda$, equals to $\rho$ identically on $\mT$. Hence $\thet$ is continuous and moreover, by Lemma \ref{lem:conv} Condition {\bf C2} from Section \ref{sec:add} holds. Therefore Proposition \ref{prop:dim=int} implies
$$
    \|\Conv_\lambda\|_\mu^2
  = \frac1{2\pi} \int_\mT \rho\, dx
  = \rho.
$$

\section{Operators of diagonal type}
\label{sec:dt}

\subsection{Definition and properties}

In this section we continue the study of bounded operators on $L^2(\mT)$. Consider the operator
$W = (W_{j,k})_{j,k\in\mZ}$
$$
  f\mapsto F = Wf, \quad
  F_j = \sum_{k\in\mZ} W_{j,k} f_k, \qquad
  f = \sum_{k\in\mZ} f_k e^{ikx}, \quad
  F = \sum_{k\in\mZ} F_k e^{ikx} .
$$

Let $\Lambda_k$, $k\in\mZ$ be the distributions
\begin{equation}
\label{Lambda_k}
  \Lambda_k(x) = \sum_{j\in\mZ} W_{k+j,j} e^{ijx} .
\end{equation}
Then by (\ref{Cnorm})
\begin{equation}
\label{ConvLambda}
  \|\Conv_{\Lambda_k}\| = c_k, \qquad  c_k = \sup_{j\in\mZ} |W_{k+j,j}| .
\end{equation}
By using the distributions $\Lambda_k$, we obtain another form of $W$ (see Remark \ref{rem:priority}):
$$
  Wf = \sum_{k\in\mZ} e^{ikx} \Lambda_k * f .
$$

\begin{dfn}
\label{dfn:||DT}
The operator $W$ is said to be of diagonal type if\;
$\sum_{k\in\mZ} c_k = {\bc} < \infty$. Below $\DT$ denotes the space of such operators.
The sequence $\{c_k\}_{k\in\mZ}$ will be said to be the majorating sequence for $W$. We also put
$$
  \|W\|_{\DT} = \bc.
$$
\end{dfn}

If $W\in\DT$ then $|W_{j,k}|\le c_{j-k}$. Hence, elements $W_{j,k}$ which are not too close to zero are concentrated near the main diagonal of the matrix $W$. This motivates the terminology.
\medskip

{\bf Examples}. 1. For any distribution (\ref{lambda(x)}) with finite $c_\lambda$ the operator $\Conv_\lambda\in\DT$ because all $c_k$ vanish except $c_0=c_\lambda$.

2. Let $g:\mT\to\mC$ be a function with absolutely converging Fourier series. Then the multiplication operator $f\mapsto \widehat g f = gf$ is of diagonal type. Indeed,
$\widehat g_{j,k}=g_{j-k}$, where $g=\sum g_k e^{ikx}$. Hence,  $c_k=|g_k|$.

We have the obvious estimate
$$
  \|g\|_\infty \le \|\widehat g\|_{\DT}.
$$

3. If $W\in\DT$ then the conjugated operator $W^*$ is also of diagonal type with the majorating sequence $\{c_k^*\}_{k\in\mZ} = \{c_{-k}\}_{k\in\mZ}$. This follows from the equation $\sup_j |W^*_{k+j,j}| = \sup_j |W_{-k+j,j}|$. Moreover, we have the equation $\|W^*\|_{\DT} = \|W\|_{\DT}$.

4. Linear combination of operators of diagonal type is an operator of diagonal type and
$$
      \|\lambda_1 W_1 + \lambda_2 W_2\|_{\DT}
  \le |\lambda_1| \|W_1\|_{\DT} + |\lambda_2| \|W_2\|_{\DT} .
$$

\begin{lem}
\label{lem:prod_diag}
Product of two operators of diagonal type is also an operator of diagonal type and
$\|W' W''\|_{\DT} \le \|W''\|_{\DT} \|W'\|_{\DT}$.
\end{lem}

{\it Proof}. Suppose $W'$ and $W''$ are of diagonal type. Let $c'_k$ and $c''_k$ be the majorating sequences. Then
$$
    |(W'W'')_{jl}|
  = \Big| \sum_k W'_{jk} W''_{kl} \Big|
 \le \tilde c_{j-l}, \qquad
     \tilde c_{j-l} = \sum_k c'_{j-k} c''_{k-l}.
$$
In particular, the majorating sequence $c_k$ of the operator $W' W''$ satisfies the inequality $c_k\le\tilde c_k$ for any $k\in\mZ$.

The obvious computation $\sum\tilde c_s = \sum_{s,k} c'_{s-k} c''_k = \bc' \bc''$ finishes the proof. \qed

\begin{lem}
\label{lem:DTclosed}
The space $\DT$ is closed w.r.t. the norm $\|\cdot\|_{\DT}$.
\end{lem}

{\it Proof}. Suppose $\{W_n\}_{n\in\mN}$ is a Cauchy sequence: for any $\eps>0$ there exists $N\in\mN$ such that
\begin{equation}
\label{implies}
  m,n>N\quad \mbox{implies}
  \quad
  \|W_m - W_n\|_{\DT} < \eps.
\end{equation}
Then for any $j,k\in\mZ$ the matrix elements $\{(W_n)_{jk}\}_{n\in\mN}$ also form a Cauchy sequence. Hence, there exists a matrix $W = \lim_{n\to\infty} W_n$, where the limit is taken element-wise:
$$
  W_{jk} = \lim_{n\to\infty} (W_n)_{jk}, \qquad j,k\in\mZ.
$$
We put
$$
  d_k = \sup_{j\in\mZ} \big| (W_n)_{k+j,j} - W_{k+j,j} \big|, \quad
  d = \|W_n - W\|_{\DT}.
$$
Then $d = \sum_{k\in\mZ} d_k$. There exists $K>0$ such that
\begin{equation}
\label{d1}
  \Big| d - \sum_{k=-K}^K d_k \Big| < \eps .
\end{equation}

There exists an integer sequence $\{j_k\}$ such that putting
$\delta_k = \big| (W_n)_{k+j_k,j_k} - W_{k+j_k,j_k} \big|$, we have:
$$
    | d_k - \delta_k |
  < \frac\eps{2K+1} \qquad
    \mbox{for any $k\in [-K,K]$},
$$
Then
\begin{equation}
\label{d2}
  \sum_{k=-K}^K |d_k - \delta_k| < \eps .
\end{equation}

Putting $\widetilde\delta_k = \big| (W_m)_{k+j_k,j_k} - (W_n)_{k+j_k,j_k} \big|$, we take $m>n$ such that
$$
  \widetilde\delta_k < \frac\eps{2K+1}, \qquad
  -K\le k\le K.
$$
Hence,
\begin{equation}
\label{d3}
  \sum_{k=-K}^K \widetilde\delta_k < \eps .
\end{equation}

By (\ref{implies})
\begin{equation}
\label{d4}
     \sum_{k=-K}^K  \big| (W_n)_{k+j_k,j_k} - (W_m)_{k+j_k,j_k} \big|
 \le \|W_n - W_m\|_{\DT}
  <  \eps .
\end{equation}
Combining (\ref{d1})--(\ref{d4}), we obtain $d\le 4\eps$. \qed

\begin{cor}
\label{cor:C*}
The space $\DT$ endowed with the norm $\|\cdot\|_{\DT}$ is a $\mC^*$-algebra.
\end{cor}

\begin{lem}
\label{lem:normdiagtype}
Suppose $W\in\DT$. Then $\|W\|\le \|W\|_{\DT}$.
\end{lem}

{\it Proof}. For any $f\in L^2(\mT)$ we have:
\begin{eqnarray}
\nonumber
     \|Wf\|^2
 &=& \sum_k \Big| \sum_j W_{kj} f_j \Big|^2
  =  \sum_{k,j,l} W_{kj} \overline W_{kl} f_j \overline f_l
 \le  \sum_{k,j,l} | W_{kj} W_{kl} f_j f_l | \\
\label{comp}
&\le& \sum_{k,j,l} c_{k-j} c_{k-l} | f_j f_l |
  =   \sum_{j,l} A_{j-l} | f_j f_l |,
\end{eqnarray}
where $A_{j-l} = \sum_{k\in\mZ} c_{k-j} c_{k-l}$. Note that
$
  \sum_{s\in\mZ} A_s = \sum_{s,k} c_k c_{s+k} = \bc^2.
$
Continuing computation (\ref{comp}), we obtain:
$$
  \|Wf\|^2 \le \sum_s A_s \sum_l |f_{l+s}|\, |f_l|
           \le \|f\|^2\sum_s A_s
           \le \bc^2 \|f\|^2.
$$
The middle inequality follows from Cauchy-Bunyakovsky-Schwarz. This implies the lemma. \qed

For any $a\in\mT$ consider the distribution $L_a$,
\begin{equation}
\label{La}
  L_a(x) = \sum_{k\in\mZ} e^{ik(x+a)} \Lambda_k(x).
\end{equation}
Then
\begin{equation}
\label{|Conv|}
      \|\Conv_{L_a}\|
  \le \sum \|\Conv_{\Lambda_k}\| = \bc.
\end{equation}
By (\ref{Lambda_k})
\begin{equation}
\label{w(a)}
     L_a(x)
  =  \sum_{j,k\in\mZ} e^{ik(x+a)} W_{k+j,j} e^{ijx}
  =  \sum_{l\in\mZ} w_l(a) e^{ilx}, \qquad
     w_l(a)
  =  \sum_{j\in\mZ} W_{l,j} e^{i(l-j)a} .
\end{equation}
For any $l\in\mZ$ and $a\in\mT$
\begin{equation}
\label{|w|}
  |w_l(a)| \le \sum_{j\in\mZ} |W_{l,j}| \le \bc.
\end{equation}

\subsection{Computation of the $\mu$-norm}

Consider the function $a\mapsto\rho(L_a)$ determined by (\ref{rho}):
\begin{equation}
\label{rho(L)}
  \rho(L_a) = \limsup_{\# I\to\infty} \rho_I(L_a), \qquad
  \rho_I(L_a) = \frac1{\# I} \sum_{l\in I} |w_l(a)|^2,
\end{equation}
where $I\subset\mZ$ are intervals.
\medskip

\begin{lem}
\label{lem:rhoLa}
If $W\in\DT$ then $\rho(L_a)$ is a continuous function of $a\in\mT$.
\end{lem}

We prove Lemma \ref{lem:rhoLa} in Section \ref{sec:La}.

\begin{prop}
\label{prop:diagonal}
Suppose $W\in\DT$. Then
\begin{equation}
\label{dimdiag}
  \|W\|_\mu^2 = \frac1{2\pi} \int_\mT \rho(L_a)\, da.
\end{equation}
\end{prop}

\begin{dfn}
We say that the matrix $(W_{k,j})$ is $\tau$-periodic, $\tau\in\mN$, if
\begin{equation}
\label{periodic}
  W_{k+\tau,j+\tau} = W_{k,j} \quad\mbox{for any } k,j\in\mZ.
\end{equation}
\end{dfn}

\begin{lem}
Let $W'$ and $W''$ be two bounded operators with $\tau$-periodic matrices. Then matrices of the operators
$$
  W^{\prime *}, \quad W'W'', \quad
  \mbox{and } a'W'+a''W'', \qquad a',a''\in\mC
$$
are also $\tau$-periodic.
\end{lem}

We skip an obvious proof.

\begin{cor}
The space of $\DT$ operators with $\tau$-periodiс matrices form a $\mC^*$ subalgebra in $\DT$.
\end{cor}

\begin{exa}
\label{exa:Wperiodic}
Suppose the matrix $(W_{k,j})$ is $\tau$-periodic, $\tau\in\mN$.
Then $w_{l+\tau}(a) = w_l(a)$ and therefore,
$\rho(L_a) = \frac1\tau \sum_{l=0}^{\tau-1} |w_l(a)|^2$.

In particular, if $\tau=1$ then $w_l=w_0$ for any $l\in\mZ$, and $W$ is the operator of multiplication by $w_0(x)$. In this case by (\ref{dimdiag})
$$
  \|W\|_\mu^2 = \frac1{2\pi} \int_\mT |w_0(a)|^2\, da = \|w_0\|^2.
$$
This result is a particular case of Lemma \ref{lem:W_g} under additional assumptions that $(\calX,\mu) = (\mT,\frac1{2\pi}dx)$ and the Fourier series of the function $g$ absolutely converges.
\end{exa}

{\it Proof of Proposition \ref{prop:diagonal}}.

\begin{lem}
\label{lem:Wf-L*f}
Suppose $W\in\DT$, $Y = [a-\eps,a+\eps]\subset\mT$ and $f=\pi_Y f\in L^2(\mT)$.
Given $\sigma>0$ for any $a\in\mT$ and any $0 < \eps \le \eps_0(\sigma,W)$
\begin{equation}
\label{Wf-L*f}
  \|Wf - L_a * f\| \le 3\sigma \|f\|.
\end{equation}
\end{lem}

{\it Proof of Lemma \ref{lem:Wf-L*f}}. We define the cut-off $W_K$,
$$
  W_K f = \sum_{|k|\le K} e^{ikx} \Lambda_k * f.
$$
Then if $K = K(\sigma)$ is large enough then
\begin{equation}
\label{cutoff}
     \Big\| Wf - W_K f \Big\|
  =  \Big\| \sum_{|k| > K} e^{ikx} \Lambda_k * f \Big\|
 \le \sum_{|k| > K} \|\Lambda_k * f \|
 \le \sum_{|k|>K} c_k \|f\|
 \le \sigma \|f\| .
\end{equation}

By (\ref{d*e})
\begin{eqnarray*}
       W_K f
  &=&  \sum_{|k|\le K} (e^{ik(x+a)}\Lambda_k) * (e^{ik(x-a)} f) = L_{a,K} * f + \Delta, \\
       L_{a,K}(x)
  &=&  \sum_{|k|\le K} e^{ik(x+a)} \Lambda_k(x), \quad
       \Delta
 \;=\; \sum_{|k|\le K} (e^{ik(x+a)} \Lambda_k) * \big((e^{ik(x-a)} - 1) f\big) .
\end{eqnarray*}
By (\ref{d*e}) and (\ref{ConvLambda})
$$
      \|\Delta\|
  \le \sum_{|k|\le K} c_k \| (e^{ik(x-a)} - 1) f \|
  \le \sum_{|k|\le K}  c_k \eps K \|f\|
   <  \eps K \bc \|f\| .
$$
Therefore
\begin{equation}
\label{WKf-L*f}
       \| W_K f - L_{a,K} * f \| \le \eps K \bc \|f\| .
\end{equation}

We also have:
\begin{equation}
\label{L-L}
     \|L_{a,K} * f - L_a * f\|
 \le \|\Conv_{L_{a,K} - L_a}\| \cdot \|f\|
 \le \sum_{|k| > K} c_k \|f\|
 \le \sigma \|f\| .
\end{equation}

Combining estimates (\ref{cutoff}), (\ref{WKf-L*f}), and (\ref{L-L}), we obtain (\ref{Wf-L*f}) if $\eps_0=\sigma/(K\bc)$. \qed

Now we return to the proof of Proposition \ref{prop:diagonal}. We associate with $W$ the function $\thet$ by (\ref{theta}). According to Lemmas \ref{lem:conv} and \ref{lem:Wf-L*f} for any $a\in\mT$ we have: $\thet(a) = \rho(L_a)$. The function $\thet$ is continuous (Lemma \ref{lem:rhoLa}) and Condition {\bf C2} from Section \ref{sec:add} holds because by Lemma \ref{lem:conv} it holds for convolutions. Hence (\ref{dimdiag}) follows from Proposition \ref{prop:dim=int}. \qed

\section{$\mu$-norm and average trace}

\subsection{Definition of $\bT(W)$}

By Lemma \ref{lem:finite} in the finite-dimensional case the $\mu$-norm coincides with the average trace of $W^* W$. The following lemma compares $\mu$-norm with the average trace for operators of diagonal type.

For any operator $W\in\DT$ we define average trace of $W^* W$ by
\begin{equation}
\label{|W|=T(W)}
    \bT(W)
  = \limsup_{\# I\to\infty} \bT(I,W), \qquad
    \bT(I,W)
  = \frac1{\# I} \sum_{j\in\mZ,\, l\in I} |W_{l,j}|^2 .
\end{equation}

\begin{prop}
\label{prop:dimDT>}
Suppose $W\in\DT$. Then
\begin{equation}
\label{bL}
    \|W\|_\mu^2 \ge \bT(W) .
\end{equation}
\end{prop}

\begin{cor}
\label{cor:.<.<.}
Combining (\ref{dim<inttheta}) and (\ref{bL}), we have:
$$
  \bT(W) \le \| W \|_\mu^2 \le \int_\mT \thet\, d\mu .
$$
\end{cor}

{\it Proof of Proposition \ref{prop:dimDT>}}. By (\ref{w(a)})
\begin{equation}
\label{w^2}
    |w_l(a)|^2
  = \sum_{j,k\in\mZ} W_{l,j} \overline W_{l,k} e^{i(k-j)a}.
\end{equation}
Then by (\ref{rho(L)}), (\ref{dimdiag})
\begin{equation}
\label{dim=intlim}
     \|W\|_\mu^2
  =  \frac1{2\pi} \int_\mT \limsup_{\# I\to\infty} \frac1{\# I}
           \sum_{j,k\in\mZ,\, l\in I} W_{l,j} \overline W_{l,k} e^{i(k-j)a} \, da
\end{equation}

Given an arbitrary big $M>0$ and arbitrarily small $\sigma>0$ let $J\subset\mZ$ be an interval such that $\# J > M$ and $\bT(J,W) > \bT(W) - \sigma$. Then by (\ref{dim=intlim})
$$
      \|W\|_\mu^2
  \ge \frac1{2\pi} \int_\mT \frac1{\# J}
             \sum_{j,k\in\mZ,\, l\in J} W_{l,j} \overline W_{l,k} e^{i(k-j)a} \, da
   =  \bT(J,W).
$$
Hence $\|W\|_\mu^2 \ge \bT(W) - \sigma$. This implies (\ref{bL}).  \qed

\subsection{Multiplication by a unitary operator}

\begin{prop}
\label{prop:T(UWU)}
Suppose $W,U\in\DT$ and $U$ is unitary. Then
\begin{equation}
\label{WU=W=UW}
  \bT(WU) = \bT(W) = \bT(UW).
\end{equation}
\end{prop}

{\it Proof}. First, note that $W' = WU$ and $W'' = UW$ are operators of diagonal type and
$$
  W'_{l,j} = \sum_{m\in\mZ} W_{l,m} U_{m,j}, \quad
  W''_{l,j} = \sum_{m\in\mZ} U_{l,m} W_{m,j}.
$$
Then the equation $U^{-1} = U^*$ implies
$$
    \sum_{j\in\mZ,\, l\in I} |W'_{l,j}|^2
  = \sum_{m,n,j\in\mZ,\, l\in I} W_{l,m} U_{m,j} \overline W_{l,n} \overline U_{n,j}
  = \sum_{m\in\mZ,\, l\in I} |W_{l,m}|^2 .
$$
Hence, the first equation (\ref{WU=W=UW}) follows.

The second equation (\ref{WU=W=UW}) requires a larger effort. Let
$U_{jl}^{-1} = \overline U_{lj}$ denote elements of the matrix $U^{-1}$. Then
$$
     \sum_{j\in\mZ} |W''_{lj}|^2
  =  \sum_{m,n,j\in\mZ} U_{l,m} W_{m,j} \overline U_{l,n} \overline W_{n,j}
  =  \sum_{m,n,j\in\mZ} U_{l,m} U_{n,l}^{-1} W_{m,j} \overline W_{n,j} .
$$

Let $I\subset\mZ$ be an interval. Since
$$
     \bT(I,W)
  =  \frac1{\# I} \sum_{l\in I,\,m,n,j\in\mZ}
         \delta_{n,l} \delta_{l,m} W_{m,jn} \overline W_{n,j},
$$
we have:
\begin{eqnarray*}
       |\bT(I,W'') - \bT(I,W)|
  &=&  \frac1{\# I} \bigg|
         \sum_{l\in I,\,m,n,j\in\mZ}
            (U_{n,l}^{-1}U_{lm} - \delta_{n,l}\delta_{lm})
                    W_{m,j}\overline W_{n,j} \bigg| \\
 &\le& \sum_{m,n,j\in\mZ} \bigg| \frac1{\# I} \sum_{l\in I}
             (U_{n,l}^{-1}U_{l,m} - \delta_{n,l}\delta_{l,m}) \bigg|
               c_{m-l} c_{n-l}
 \; =\; \Delta,  \\
        \Delta
   =   \sum_{m,n\in\mZ} \Gamma_{m,n} \widetilde c_{m-n},
  \!\!\!\!\!\!\!\!
  &&  \quad
       \Gamma_{m,n}
   =   \bigg| \frac1{\# I} \sum_{l\in I}
             (U_{n,l}^{-1}U_{l,m} - \delta_{n,l}\delta_{l,m}) \bigg|,
\end{eqnarray*}
where $\widetilde c_{m-n} = \sum_{j\in\mZ} c_{m-l} c_{n-l}$,\;
$\sum_{m\in\mZ} \widetilde c_m = \bc^2$.

Since $U,U^{-1}\in\DT$, we have:
$|U_{kl}| \le d_{l-k}$\;, $\sum_{j\in\mZ} d_j = \bd < \infty$.
Given $\sigma>0$ there exists $M>0$ such that
$$
  \sum_{|j|>M} d_j < \sigma.
$$
We put $\displaystyle e_{n-l} = \sum_{m\in\mZ} d_{m-l} \widetilde c_{m-n}$,
$$
  I_M^+ = \{k : [k-M,k+M]\cap I \ne \emptyset\}, \quad
  I_M^- = \{k : [k-M,k+M] \subset I\} .
$$
Then $\sum_{j\in\mZ} e_j = \bd\bc^2$.

We have: $\Delta = \Delta_1 + \Delta_2 + \Delta_3$:
\begin{eqnarray*}
      \Delta_1
  &=& \sum_{m\in\mZ,\,n\not\in I^+_M} \Gamma_{m,n}
               \widetilde c_{m-n}, \\
      \Delta_2
  &=& \sum_{m\in\mZ,\,n\in I^-_M} \Gamma_{m,n}
               \widetilde c_{m-n}, \\
      \Delta_3
  &=& \sum_{m\in\mZ,\,n\in I^+_M\setminus I^-_M}
           \Gamma_{m,n}
               \widetilde c_{m-n} .
\end{eqnarray*}

First, we estimate $\Delta_1$:
\begin{eqnarray*}
       \Delta_1
 &\le& \frac1{\# I} \sum_{m\in\mZ,\,n\not\in I^+_M} \sum_{l\in I}
             d_{n-l} d_{m-l} \widetilde c_{m-n} \\
 &\le& \frac1{\# I} \sum_{n\not\in I^+_M,\,l\in I}
             d_{n-l} e_{n-l}
 \; \le \; \bd\bc^2 \sigma .
\end{eqnarray*}

Now we turn to $\Delta_2$. Since
$\sum_{l\in\mZ} U^{-1}_{n,l} U_{l,m}
 = \sum_{l\in\mZ} \delta_{n,l}\delta_{m,l} = \delta_{n,m}$, we have:
\begin{eqnarray*}
      \Delta_2
  &=& \sum_{m\in\mZ,\,n\in I^-_M} \bigg| \frac1{\# I} \sum_{l\in I}
        (U^{-1}_{n,l}U_{l,m} - \delta_{n,l}\delta_{m,l}) \bigg|\,\widetilde c_{m-n}, \\
  &=& \frac1{\# I} \sum_{m\in\mZ,\,n\in I^-_M}
        \bigg| \sum_{l\in\mZ\setminus I} U^{-1}_{n,l}U_{l,m} \bigg|\,
                   \widetilde c_{m-n}, \\
 &\le& \frac1{\# I} \sum_{m\in\mZ,\,n\in I^-_M,\,l\in\mZ\setminus I}
            d_{n-l} d_{m-l} \widetilde c_{m-n}, \\
  &=& \frac1{\# I} \sum_{n\in I^-_M,\,l\in\mZ\setminus I}
            d_{n-l} e_{n-l}
 \; \le \; \sum_{|m| > M} d_m e_m
 \; \le \; \bd\bc^2 \sigma.
\end{eqnarray*}

Finally, we estimate $\Delta_3$:
\begin{eqnarray*}
       \Delta_3
 &\le& \frac1{\# I} \sum_{m\in\mZ,\,n\in I^+_M\setminus I^-_M,\,l\in I}
             (d_{n-l} d_{m-l} + \delta_{n,l}\delta_{l,m})
               \widetilde c_{m-n} \\
 &\le& \frac1{\# I} \sum_{n\in I^+_M\setminus I^-_M,\,l\in I}
             (d_{n-l} e_{n-l} + \delta_{n,l} \widetilde c_{l-n}) \\
  &=&  \frac1{\# I} \sum_{n\in I^+_M\setminus I^-_M}
             (\bd^2\bc^2 + \bc^2)
 \; =\; \frac{2(M+1)}{\# I} (\bd^2 + 1) \bc^2.
\end{eqnarray*}

We see that $\Delta$ can be made arbitrarily small if we choose sufficiently small $\sigma$ and $M / \# I$. This implies
$$
    \limsup_{\# I\to\infty} \bT(I,W)
  = \limsup_{\# I\to\infty} \bT(I,W'').
$$
\qed

\section{Technical statements}

\subsection{Norm of a function on a subset}

\begin{lem}
\label{lem:piDf}
Let $D\subset A\subset\calX$ be two measurable sets and let $f = \pi_A f \in L^2(\calX)$ satisfy $c<\big|f|_A\big|<c^{-1}$. Then
\begin{equation}
\label{piDf}
       \|\pi_D f\|
  \le  \frac{\sqrt{\mu(D)}}{c^2 \sqrt{\mu(A)}} \|f\| .
\end{equation}
\end{lem}

{\it Proof}. We have: $\|f\|^2 \ge c^2\mu(A)$ and $\|\pi_D f\|^2 \le c^{-2}\mu(D)$. Therefore
$$
  \|\pi_D f\|^2 \le c^{-2}\mu(D) \frac{\|f\|^2}{c^2 \mu(A)}.
$$
This implies (\ref{piDf}). \qed

\subsection{The function $\rho(\lambda)$}

In this section we present two lemmas on the function $\rho=\rho(\lambda)$ for a distribution $\lambda$, determined by (\ref{lambda(x)}).

\begin{lem}
\label{lem:C<J<2C}
For any $\beta>0$ and any $C>1$ there exists an interval $J\subset\mZ$ such that
$$
  \rho_J(\lambda) > \rho - \beta \quad
  \mbox{and}\quad
  C \le \# J \le 2C + 1.
$$
\end{lem}

{\it Proof}. There exists $J_0$ such that
$$
  \rho_{J_0}(\lambda) > \rho - \beta \quad
  \mbox{and}\quad
  \# J_0 > 2C + 1.
$$
We break $J_0$ into two nonintersecting intervals $J'_0$ and $J''_0$ such that $|\# J'_0 - \# J''_0| \le 1$. One of these intervals (we redenote it by $J_1$) satisfies $\rho_{J_1}(\lambda) > \rho - \beta$. If
$\# J_1 > 2C + 1$, we repeat the argument. \qed

\begin{lem}
\label{lem:rhoJ0}
Let $M'=M'(\sigma)$ satisfy (\ref{M'}). Let $\widetilde J$ be an integer interval such that
\begin{equation}
\label{rhohatJ}
  \rho_{\widetilde J}(\lambda) \ge \rho - \widetilde\sigma, \quad
  \#\widetilde J > 2M'.
\end{equation}
Let $J\subset\widetilde J$ be the interval which is obtained from $\widetilde J$ by removing two intervals: extreme left $J_l$ and extreme right $J_r$, $\# J_l = \# J_r = M'$.

Then for any integer interval $J_0\subset J$
\begin{equation}
\label{rhoJ0}
  \rho_{J_0}(\lambda) = \rho - \sigma_0, \quad
  \sigma_0 \le (\widetilde\sigma + \sigma) \frac{\#\widetilde J}{\# J_0} - \sigma.
\end{equation}
\end{lem}

{\it Proof}. We have: $\widetilde J = \widetilde J_l \cup J_0\cup \widetilde J_r$, where $\widetilde J_l\supset J_l$ and $\widetilde J_r\supset J_r$ are integer intervals and $\#\widetilde J_l\ge M'$,
$\#\widetilde J_r\ge M'$. Hence, by (\ref{rhohatJ})
$$
      \#\widetilde J (\rho - \widetilde\sigma)
  \le \sum_{j\in\widetilde J} |\lambda_j|^2
  \le (\rho - \sigma_0)\# J_0 + (\rho + \sigma) (\#\widetilde J_l + \#\widetilde J_r).
$$
Therefore
$\# J_0\sigma_0 \le \#\widetilde J\widetilde\sigma + (\#\widetilde J - \# J_0)\sigma$. This implies (\ref{rhoJ0}). \qed

\subsection{The function $\rho(L_a)$}
\label{sec:La}

In this section we study the distribution $L_a$, determined by (\ref{La}).

\begin{lem}
\label{lem:equicont}
Suppose $W\in\DT$. Then the functions $a\mapsto\rho_I(L_a)$ are equicontinuous i.e., for any $\sigma>0$ there exists $\gamma(\sigma)$ such that
$$
  |\rho_I(L_a) - \rho_I(L_b)| < \sigma \quad
  \mbox{whenever $|a-b|<\gamma(\sigma)$}.
$$
\end{lem}

It is important that $\gamma$ does not depend on $I$.

{\it Proof of Lemma \ref{lem:equicont}}. By (\ref{w(a)}) and (\ref{rho(L)}) we have the estimate
\begin{eqnarray*}
     |\rho_I(L_a) - \rho_I(L_b)|
 &=& \frac1{\# I} \Big| \sum_{l\in I} \sum_{j,n\in\mZ} W_{l,j} \overline W_{l,n}
                                (e^{i(n-j)a} - e^{i(n-j)b}) \Big| \\
&\le& \frac1{\# I} \sum_{l\in I,\,j,n\in\mZ} c_{l-j} c_{l-n}
                               | (e^{i(n-j)a} - e^{i(n-j)b} |  \\
 &=&  \sum_{j,n\in\mZ} c_{n} c_j | (e^{i(n-j)a} - e^{i(n-j)b} |
\;\le\; \Sigma_1 + \Sigma_2,
\end{eqnarray*}
where
$$
  \Sigma_1 = \sum_{|j-n|\le K} c_n c_j | (e^{i(n-j)a} - e^{i(n-j)b} | , \quad
  \Sigma_2 = \sum_{|j-n| > K} c_n c_j | (e^{i(n-j)a} - e^{i(n-j)b} | .
$$

We have:
\begin{eqnarray*}
       \Sigma_1
 &\le& K|a-b| \sum_{|j-n|\le K} c_n c_j
  \le  \bc^2 K |a-b|, \\
       \Sigma_2
 &\le& 2 \sum_{|j-n| > K} c_n c_j
  \le  2 \sum_{j\in\mZ, |n| > K/2} c_n c_j + 2 \sum_{n\in\mZ, |j| > K/2} c_n c_j
  \le  4\bc \sigma_0(K), \\
       \sigma_0(K)
  &=&  \sum_{|n| > K/2} c_n .
\end{eqnarray*}
Hence, for $|a-b|<\gamma$
$$
      |\rho_I(L_a) - \rho_I(L_b)|
  \le \bc^2 K\gamma + 4\bc\sigma_0(K).
$$
To make the right-hand side of this inequality arbitrarily small, we first, choose sufficiently large $K$ and then take $\gamma=\gamma(K)>0$ sufficiently small.  \qed

Now Lemma \ref{lem:rhoLa} turns out to be a simple corollary of Lemma \ref{lem:equicont}. For completeness we present the argument.

Given $\sigma>0$ take an interval $I\subset\mZ$ with arbitrarily large $\# I$ and such that $\rho_I(L_a)>\rho(L_a)-\sigma$. Then for any $b\in\mT$ such that $|b-a|<\gamma(\sigma)$, where $\gamma(\sigma)$ is defined in Lemma \ref{lem:equicont}, we have:
$$
  \rho_I(L_b) \ge \rho_I(L_a) - \sigma > \rho(L_a) - 2\sigma.
$$

Analogously if $N(\sigma)$ is sufficiently large then for any $I\subset\mZ$ with $\# I > N(\sigma)$ we have: $\rho_I(L_a)<\rho(L_a)+\sigma$. Then for any $b\in\mT$ satisfying $|b-a|<\gamma(\sigma)$
$$
  \rho_I(L_b) \le \rho_I(L_a) +\sigma < \rho(L_a) + 2\sigma.
$$
Hence $|\rho(L_a) - \rho(L_b)| < 2\sigma$ whenever $|a-b| < \gamma(\sigma)$.

\end{document}